\documentclass[11pt]{article}

\usepackage[top=1in, bottom=1in, left=1in, right=1in]{geometry}

\usepackage{amssymb,graphics,amsmath,amsthm,amsopn,amstext,amsfonts,color,fullpage,fancybox,hyperref,bm}
\usepackage{graphicx}
\usepackage{verbatim}
\usepackage{algorithm}
\usepackage{algpseudocode}
\usepackage{longtable}
\usepackage{multirow}
\usepackage{diagbox}
\usepackage{lscape}

\newtheorem{assumption}{Assumption}
\newtheorem{remark}{Remark}
\newtheorem{theorem}{Theorem}
\newtheorem{lemma}{Lemma}
\newtheorem{corollary}{Corollary}
\newtheorem{proposition}{Proposition}
\newtheorem{definition}{Definition}

\newcommand{\epc}{\hspace{1pc}}
\def\norm#1{\big\| #1\big\|}
\def\inprod#1#2{\big\langle #1,\,#2\big\rangle}

\newcommand{\bea}{\begin{eqnarray*}}
	\newcommand{\eea}{\end{eqnarray*}}

\def\sig{\sigma}

\def\inprod#1#2{\langle#1,\,#2\rangle}
\def\norm#1{\|#1\|}

\def\grad{\nabla}

\def\S{{\cal S}}

\def\S{\mathbb{S}}

\def\Res{{\cal R}}

\def\R{\mathbb{R}} 
\def\S{\mathbb{S}}

\def\X{\mathbb{X}}
\def\Y{\mathbb{Y}}

\def\W{\mathbb{W}}

\def\EMAIL#1{\href{mailto:#1}{#1}}

\begin{document}

\title{\Large On Degenerate Doubly Nonnegative Projection Problems}

\author{
Ying Cui \thanks{Department of Industrial and Systems Engineering, University of Minnesota, Minneapolis, USA. \EMAIL{yingcui@umn.edu}}, \;
Ling Liang \thanks{Department of Mathematics, National University of Singapore, Singapore.\EMAIL{liang.ling@u.nus.edu}},\;
Defeng Sun \thanks{Department of Applied Mathematics, The Hong Kong Polytechnic University, Hung Hom, Hong Kong. \linebreak \EMAIL{defeng.sun@polyu.edu.hk}}, \; and\;
Kim-Chuan Toh, \thanks{Department of Mathematics, and Institute of Operations Research and Analytics, National University of Singapore, Singapore. \EMAIL{mattohkc@nus.edu.sg}}}

\date{}

\maketitle

\begin{abstract}
The doubly nonnegative (DNN) cone, being the set of all positive semidefinite matrices whose elements are nonnegative, is a popular approximation of the computationally intractable completely positive cone. The major difficulty for implementing a Newton-type method to compute the projection of a given large scale matrix onto the DNN cone lies in the possible failure of the constraint nondegeneracy, a generalization of the linear independence constraint qualification for nonlinear programming. Such a failure results in the singularity of the Jacobian of the nonsmooth equation representing the Karush-Kuhn-Tucker optimality condition  that prevents the semismooth Newton-CG method from solving it with a desirable convergence rate. In this paper,  we overcome the aforementioned difficulty by solving a sequence of better conditioned nonsmooth equations generated by the augmented Lagrangian method (ALM) instead of solving  one above mentioned singular equation. By leveraging on the metric subregularity of the normal cone associated with the positive semidefinite cone, we derive sufficient conditions to ensure the dual quadratic growth condition of the underlying problem, which further leads to the asymptotically superlinear convergence  of the proposed ALM.  Numerical results on difficult randomly generated instances and from the semidefinite programming library  are presented to demonstrate the efficiency of the algorithm for computing the DNN projection to a very high accuracy.
\end{abstract}

{\small
\begin{center}
\parbox{0.95\hsize}{{\bf Keywords.}\;  Doubly nonnegative cone, semidefinite programming, augmented Lagrangian method, semismooth Newton, degeneracy, metric subregularity}
\end{center}
}
\begin{center}
\parbox{0.95\hsize}{{\bf AMS subject classifications:}\; 90C06, 90C22, 90C25}
\end{center}

\maketitle

\section{Introduction}\label{intro}
Let $\mathbb{S}^n$ be the vector space of $n\times n$ symmetric matrices,
$\mathbb{S}_+^n$ the cone of $n\times n$ symmetric positive semidefinite matrices, $\mathbb{R}_+^n$ the nonnegative orthant in $\mathbb{R}^n$, and
$\mathbb{N}^n$ the nonnegative orthant in $\mathbb{R}^{n\times n}$.
The cone of $n\times n$ copositive matrices, and its dual cone, the cone of $n\times n$ completely positive matrices, are given respectively by
$$
\mathbb{C}^n\,\triangleq \,\{\,X\in\mathbb{S}^n\mid a^{\,T}Xa\geq 0, \; \forall\; a\in\mathbb{R}_+^n\,\}\quad \mbox{and}\quad
\mathbb{C}^{n,*}\,\triangleq \, {\rm conv}\,\{ \, aa^{\,T} \mid a\in\mathbb{R}_+^n\, \},
$$
where $\mbox{conv}\,\{\,C\}$ denotes the convex hull of a given set $C$.
Copositive and completely positive cones have received considerable attentions in recent years as many combinatorial and nonconvex quadratic optimization problems can be formulated equivalently as linear conic programming problems over $\mathbb{C}^n$ or $\mathbb{C}^{n,*}$, see, e.g., \cite{Bomze2002on,Burer2009,deklerk2002approximation,Povh2007copositive,Povh2009copositive}. However, both cones are computationally intractable, in the sense that  to check whether a given matrix lies in $\mathbb{C}^n$ is  co-NP-complete \cite{Murty1987}  and in $\mathbb{C}^{n,*}$ is NP-hard \cite{DickinsonGijben14}. One may refer to the survey paper \cite{Dur2010} for further properties of these two cones. A popular relaxation of the completely positive cone is  the following  doubly nonnegative  (DNN) cone
$$
\mathbb{D}^n\,\triangleq \,\{X\in\mathbb{S}^n\mid X\in\mathbb{S}_+^n, \; X\in\mathbb{N}^n\}.
$$
Clearly we have $ \mathbb{C}^{n,*}\subseteq \mathbb{D}^n$. The equality in this relation holds for $n\leq 4$ and the inclusion is strict if $n\geq 5$~\cite{MM1962}.

In this paper, we focus on designing an efficient solver for computing the projection of a given matrix $G\in\mathbb{S}^n$ over the DNN cone, i.e., finding the optimal solution of the following convex optimization problem
\begin{equation}\label{eq:DNN}
\tag{\bf P}
\operatornamewithlimits{minimize}_{X\in \mathbb{S}^n}\; \left\{\, \frac{1}{2}\|X - G\|^2 \;\left|\right. \; X\in\mathbb{S}_+^n,\; X\in \mathbb{N}^n\, \right\},
\end{equation}
where $\norm{\cdot}$ denotes the Frobenius norm. As a basic building block of various algorithms for DNN conic programming problems, such as the one in \cite{Kim2016Lagrangian}, the efficient computation of the projection onto a DNN cone is an important problem of considerable interest. For example, an efficient routine for computing this projector can be embedded in the projected gradient method for solving
$$
\operatornamewithlimits{minimize}_{X\in \mathbb{S}^n} \; \left\{\, f(X) \, \left|\right.\, X\in \mathbb{D}^n\, \right\},
$$
with $f$ being a possibly nonsmooth nonconvex function.

In a series of works \cite{li2015qsdpnal,yang2015sdpnal+,zhao2010newton}, a semismooth Newton-CG based dual augmented Lagrangian method (ALM) is proposed to solve the class of linear and convex quadratic semidefinite programming (SDP) problems. The algorithm performs fairly well for large-scale nondegenerate (note that in this paper  the concept of degeneracy refers to the constraint degeneracy for optimization problems, see for instance Definition \ref{defn:cq}) SDP problems with the dimension of the matrix variable $n$ being in the range of a few thousands but the number of equality constraints can be in the range of a few millions. If a large number of linear inequality constraints (such as the entrywise nonnegativity of the variables) are also added to the linear and convex quadratic SDP problems, it is highly possible that multiple dual solutions exist such that the generalized Jacobians of semismooth equations corresponding to the optimality conditions of augmented Lagrangian subproblems are singular. Consequently, a semismooth Newton method applied to solved the subproblems may not have fast local convergence. To resolve this issue, a majorized ALM is employed in \cite{li2015qsdpnal,yang2015sdpnal+}, where the degenerate multi-block ALM subproblems are solved by a block coordinate descent decomposition method for which each of its steps solves a nondegenerate problem involving a single block. A similar decomposition idea is adopted in \cite{CST2018} to  compute the best approximation problem over the intersection of a polyhedral set and the DNN cone.

The degeneracy issue also happens to the DNN projection problem. The dual of \eqref{eq:DNN} takes the form of
\begin{equation}\label{dual}
\tag{\bf D}
-\operatornamewithlimits{minimize}_{S, \, Z\in \mathbb{S}^n} \; \left\{\, \frac{1}{2}\|S+Z + G\|^2 - \frac{1}{2}\|G\|^2 \;\left|\right.\; S\in \mathbb{S}_+^n\, ,\; Z\in \mathbb{N}^n\,\right\}.
\end{equation}
A notable feature of the DNN projection problem is that multiple solutions to \eqref{dual} may exist, especially when the solution to \eqref{eq:DNN} possesses both the low rank and sparse properties, making the problem \eqref{eq:DNN} constraint degenerate (see Section \ref{sec:degeneracy} for detailed discussions on this part). This feature indicates the high possibility for the singularity of the generalized Jacobian of the nonsmooth equations representing the Karush-Kuhn-Tucker (KKT) optimality condition of \eqref{eq:DNN}:
\begin{equation}\label{defn:kkt}
\Res(X,S,Z)\,\triangleq \,\left(\begin{array}{cc}
X - G - S- Z\\[3pt]
X - \Pi_{\mathbb{S}_+^n}(X-S)\\[3pt]
X - \Pi_{\mathbb{N}^n}(X-Z)\\[3pt]
\end{array}
\right) = 0,\quad X,S,Z\in \mathbb{S}^n,
\end{equation}
where $\Pi_C(\cdot)$ denotes the metric projection onto a given closed convex set $C$.
It is known that the convergence rate of the conjugate gradient (CG) method for solving a linear equation is determined by its condition number. Therefore, even though the above equation is semismooth \cite{Sun2002semismooth}, directly solving it by the semismooth Newton-CG method seems not suitable when degeneracy occurs.

An important property of \eqref{eq:DNN} that distinguishes it from general convex quadratic SDP problems is the strong convexity of the objective function, which implies the uniqueness of its primal optimal solution. This motivates us to consider a primal ALM to solve the problem. Let $\sigma$ be a given positive penalty parameter. The augmented Lagrangian function of \eqref{eq:DNN} is given by
$$
L_{\sigma}(X;S,Z)\,\triangleq \,\frac{1}{2}\|X-G\|^2 + \frac{1}{2\sigma}\left(\,\|\, \Pi_{\mathbb{S}_+^n}\left(S - \sigma X\right)\, \|^2 + \|\, \Pi_{\mathbb{N}^n}\left(Z - \sigma X\right)\, \|^2\,\right), \quad X,S,Z\in\mathbb{S}^n.
$$
Given a sequence of positive scalars $\sigma_k\uparrow \sigma_\infty\leq +\infty$, the $(k+1)$-th iteration of the ALM takes the form of 
\begin{equation}\label{ALM:iter}
\left\{\begin{array}{ll}
X^{\,k+1} \approx  \displaystyle\operatornamewithlimits{argmin}_{X\in \mathbb{S}^n}\;  \left\{\,f_k(X)\,\triangleq\, L_{\sigma_k}\left(X;S^{\,k},Z^{\,k}\right)\,\right\},\\[0.15in]
(S^{\,k+1},Z^{\,k+1}) = \left(\, \Pi_{\mathbb{S}_+^n}\left(S^{\,k}  - \sigma_k X^{\,k+1}\right), \, \Pi_{\mathbb{N}^n}\left(Z^{\,k}  - \sigma_k X^{\,k+1}\right)\, \right),
\end{array}\right. \quad k\geq 0.
\end{equation}
Obviously, the major computational cost of the above framework comes from the computation of the  approximate solutions of the subproblems. The optimality condition of these subproblems can be characterized by the semismooth equations
\begin{equation}\label{eq:semismooth Newton}
0\;=\; \nabla f_k(X) \;=\; X - G -
\Pi_{\mathbb{S}_+^n}(S^k - \sigma_k X) -
\Pi_{\mathbb{N}^n}(Z^k - \sigma_k X).
\end{equation}
Different from the semismooth equation \eqref{defn:kkt}, the generalized Jacobian of the above equation is always nonsingular at any point in $\mathbb{S}^n$ (see Section \ref{sec:SSN} for the expression of its generalized Jacobian). Thus, instead of solving one singular nonsmooth KKT equation \eqref{defn:kkt}, we adopt the Newton-CG method to  solve a sequence of nonsingular nonsmooth equations \eqref{eq:semismooth Newton}.

Given the promising convergence rate of the inner semismooth Newton-CG method, the overall performance of the above proposed method depends heavily on the convergence rate of the outer augmented Lagrangian iterations. In a recent work \cite{CST2017}, it was shown that the KKT residual of the sequence of iterates generated by ALM converges asymptotically superlinearly under the dual quadratic growth condition. Though the dual quadratic growth condition has been shown to hold under  the dual second order sufficient condition \cite[Theorem 3.137]{bonnans2013perturbation}, a  unique dual optimal solution has to exist in order to fulfill the latter condition. In this paper, we show that when the dual problem has multiple solutions, the existence of a strict complementarity solution also implies that such a dual quadratic growth condition holds at any dual solution. Besides applying to this particular problem, the established theory in this paper, together with that in \cite{CST2017}, also partially explains why the ALM usually outperforms first order methods for solving other types of SDP problems to high accuracy.

In summary, the contributions of our paper are two-fold:

$\bullet$ Theoretically, we provide  sufficient conditions to ensure the quadratic growth condition of a general class of linearly constrained convex problems involving non-polyhedral functions, which includes \eqref{eq:DNN} as a special case. Besides its independent interest in variational analysis,  the derived results provide sufficient conditions for the asymptotic superlinear convergence of the KKT residual generated by the iterative sequence from the ALM.

$\bullet$ Numerically, we develop an efficient solver for computing the projection of a given matrix onto the doubly nonnegative cone to a very high accuracy. We conduct rigorous numerical experiments on various SDP instances to demonstrate the effectiveness of the proposed method.

The remaining parts of this paper are organized as follows. In the next section, we discuss necessary conditions for the constraint nondegeneracy of problem \eqref{eq:DNN} and a consequence of its failure for the Newton-type algorithm. This motivates us to  consider the Newton-CG based augmented Lagrangian method in \eqref{ALM:iter} to solve \eqref{eq:DNN}. Section \ref{sec:growth condition} is devoted to extensive studies on sufficient conditions for the quadratic growth condition of linearly constrained convex SDPs, which include \eqref{eq:DNN} as a special case. Such a quadratic growth condition ensures the asymptotically superlinear convergence rate of the proposed ALM. In Section \ref{sec:SSN}, we introduce a semismooth Newton-CG based ALM and show how it overcomes the degeneracy of the DNN constraints. Extensive numerical experiments are conducted in Section \ref{sec:numerical} to demonstrate the effectiveness of the proposed method. We conclude our paper  in the final section.

Below we list the notation to be used in our paper.
\begin{itemize}
	\item We use $\mathbb{U}$, $\mathbb{V}$, $\W$, $\X$, $\Y$ and $\mathbb{Z}$ to denote  finite dimensional real Euclidean spaces each equipped with an inner product $\langle \cdot, \cdot \rangle$ and its induced norm $\norm{\cdot}$.
	
	\item Let $\alpha \subseteq \{1,...,m\}$ and $\beta \subseteq \{1,...,n\}$ be two index sets. For any $Z \in \R^{m\times n}$, we write $Z_{\alpha\beta}$ {to be} the $|\alpha|\times |\beta|$ sub-matrix of $Z$ obtained by removing all the rows of $Z$  not in $\alpha$ and all the columns of $Z$  not in $\beta$. We denote $\textup{diag}(x_\alpha)$ as the $|\alpha|\times|\alpha|$ diagonal matrix whose diagonal entries are those of $x_{\alpha}$.
	
	\item Let $D\subseteq \X$ be a set.  For any $x\in \X$,	define ${\rm dist}(x,D)\,\triangleq \, \inf_{d\in D} \|x-d\|$. {We let}  $\delta_{D}(\cdot)$ to	be the indicator function over $D$, i.e., $\delta_{D}(x) = 0$ if $x\in D$, and $\delta_D(x) = \infty$ if $x\not\in D$. If $D\subseteq \X$ is a convex set, we  use $\textup{ri}(D)$ to denote its relative interior. For a given closed convex set $D\subseteq \X$, the metric projection of $x\in\X$ onto $D$ is defined by $\Pi_{D}(x) \,\triangleq \, \arg\min \{\|x-d\|\mid d\in D\}$. For any $x\in D$, we use $\mathcal{T}_{D}(x)$ and $\mathcal{N}_{D}(x)$ to denote the tangent and normal cones of $D$ at $x$,  and ${\rm lin}(D)$ as the lineality space of $D$, i.e., the largest linear subspace in $D$. If $D$ is a closed convex cone, we use $D^\circ$ and $D^*$ to denote the polar of $D$ and  the dual of $D$, respectively, i.e., $D^\circ\,\triangleq \,\{x\in \mathbb{X}\mid \langle x, d\rangle \leq 0,\; \forall\,d\in D\}$ and $D^*\,\triangleq \,-D^\circ$.
	
	\item For any set-valued mapping $\Gamma:\mathbb{U}\rightrightarrows\mathbb{V}$, we use $\text{gph}\,\Gamma$ to denote the graph of $\Gamma$, i.e., $\text{gph}\,\Gamma\,\triangleq \, \{(u,v)\in\mathbb{U}\times \mathbb{V}\,\mid\,v\in \Gamma(u)\}$. For any $\bar{u}\in\mathbb{U}$ and $\varepsilon>0$, denote $\mathbb{B}_{\varepsilon}(\bar{u})\,\triangleq \,\{u\in\mathbb{U}\mid \|u- \bar{u}\|\leq \varepsilon\}$.
\end{itemize}

\section{A Consequence of the  Constraint Degeneracy}\label{sec:degeneracy}

In this section, we provide necessary conditions for the primal constraint nondegeneracy and a consequence of its failure when designing Newton-type algorithms.

We start with the formal definition of the constraint nondegeneracy. Let $\mathcal{K}$ be a closed convex set in $\mathbb{Y}$. The tangent cone of $\mathcal{K}$ at  a point $y\in \mathcal{K}$ is defined by
\[
\mathcal{T}_\mathcal{K}(y) =  \left\{\,d\in \mathbb{Y} \,\mid\, \mbox{dist}\,(y+td,\mathcal{K}) = o(t),\; t\geq 0\,\right\}.
\]

Let $f: \mathbb{X} \to \mathbb{R}$ be a twice continuously differentiable function, $G: \mathbb{X} \to \mathbb{Y}$ be a twice continuously differentiable mapping and $\mathcal{K}$ be a closed convex set in $\mathbb{Y}$. For the conic programming with the form
\begin{equation}\label{gen:opt}
\operatornamewithlimits{minimize}_{x\in \mathbb{X}}\, f(x), \quad \mbox{subject to} \epc G(x)\in \mathcal{K},
\end{equation}
we have the following definition of constraint nondegeneracy \cite{Robinson03}.
\begin{definition}\label{defn:cq}
	We say that a feasible point $\bar{x}\in\mathbb{X}$ to \eqref{gen:opt} is constraint nondegenerate if
	$$
	G'(\bar{x})\,\mathbb{X} + {\rm lin}(\mathcal{T}_{\mathcal{K}}(G(\bar{x}))) = \mathbb{Y},
	$$
	where $G^\prime(\bar{x})$ denotes the Jacobian of $G$ at $\bar{x}$ and $\mbox{lin}(S)$ denotes the lineality space of a given set $S$. We say that a feasible point $\bar{x}$ is constraint degenerate if the above condition fails at $\bar{x}$.
\end{definition}

The constraint nondegeneracy condition above reduces to the linear independence constraint qualification when the problem \eqref{gen:opt} is a conventional nonlinear programming problem \cite{Robinson84, Shapiro03}. One may refer to the monograph \cite{bonnans2013perturbation} for more discussions on this concept in the context of conic program. Based on Definition \ref{defn:cq}, the constraint nondegeneracy is said to hold at a feasible point $\overline{X}\in \mathbb{S}^n$ to \eqref{eq:DNN} if
\begin{equation}\label{LICQ}
\left(
\begin{array}{cc}
\mathcal{I} \\[0.1in] \mathcal{I}
\end{array}
\right)\mathbb{S}^n + \left(
\begin{array}{cc}
{\rm lin}\left(\mathcal{T}_{\,\mathbb{S}_+^n}\left(\,\overline{X}\,\right)\right) \\[0.1in] {\rm lin}\left(\mathcal{T}_{\,\mathbb{N}^n}\left(\,\overline{X}\,\right)\right)
\end{array}
\right) =
\left(
\begin{array}{cc}
\mathbb{S}^n\\[0.1in]
\mathbb{S}^n
\end{array}\right),
\end{equation}
where $\mathcal{I}:\mathbb{S}^n\to\mathbb{S}^n$ is the identity map in $\mathbb{S}^n$.

For any given $\overline{X}\in \mathbb{S}_+^n\cap \mathbb{N}^n$, suppose that it has the following eigenvalue decomposition:
\begin{equation}\label{defn:rank}
\overline{X} = [P_{\alpha}\,P_{\bar{\alpha}}]\,{\rm diag}(\lambda_1, \lambda_2 \ldots, \lambda_r,0,\ldots, 0)\,[P_{\alpha}\,P_{\bar{\alpha}}]^{\,T},
\end{equation}
where $\alpha=\{1, 2, \ldots, r\}$, $\bar{\alpha} = \{ r+1,\ldots,n\}$,
$\lambda_1\geq \lambda_2\geq \ldots\geq \lambda_r >0$
are the positive eigenvalues of $\overline{X}$, and $P = [P_{\alpha},\,P_{\bar{\alpha}}]\in \mathcal{O}^n$ is a corresponding orthogonal matrix of orthonormal eigenvectors. We also denote
\begin{equation}\label{nnz}
\mathcal{E} = \{(i, j) \mid  \overline{X}_{ij} > 0,\; 1\leq i\leq j\leq n\}, \quad \overline{\mathcal{E}} = \{(i, j)\mid  \overline{X}_{ij} = 0,\; 1\leq i\leq j\leq n\}.
\end{equation}
It can be easily checked that (see, e.g., \cite{Arnold1971})
$$\left\{\begin{array}{ll}
\mathcal{T}_{\,\mathbb{S}_+^n}\left(\,\overline{X}\,\right)=\left\{H\in \mathbb{S}^n\mid P_{\bar\alpha}^{\,T} H\, P_{\bar\alpha}\succeq 0 \right\}, \quad \mathcal{T}_{\,\mathbb{N}^n}\left(\,\overline{X}\,\right)=\left\{H\in \mathbb{S}^n\mid H_{\overline{\mathcal{E}}}=H^{\,T}_{\overline{\mathcal{E}}}\geq 0 \right\}, \\[0.1in]
{\rm lin}\left(\mathcal{T}_{\,\mathbb{S}_+^n}\left(\,\overline{X}\,\right)\right)= \left\{H\in \mathbb{S}^n\mid P_{\bar\alpha}^{\,T} H\,P_{\bar\alpha}=0 \right\},\quad {\rm lin}\left(\mathcal{T}_{\,\mathbb{N}^n}\left(\,\overline{X}\,\right)\right)=\left\{H\in \mathbb{S}^n\mid H_{\overline{\mathcal{E}}} =H^{\,T}_{\overline{\mathcal{E}}} = 0\right\}.
\end{array}\right.
$$

The following proposition characterizes the constraint nondegeneracy of the DNN projection problem  \eqref{eq:DNN}. A necessary condition for $\overline{X}\in \mathbb{S}^n$ to be constraint nondegenerate in terms of its rank and cardinality then follows easily.

\begin{proposition}\label{prop:LICQ}
	Let $\overline{X}\in \mathbb{S}^n$ be a feasible point to \eqref{eq:DNN} with the index sets $\alpha$ and $\mathcal{E}$ given in \eqref{defn:rank} and \eqref{nnz}, respectively. Then $\overline{X}$ is constraint nondegenerate if and only if
	\begin{equation}\label{prop:licq:eqiv}
	\left\{ H \in \mathbb{S}^n \mid H_{\mathcal{E}} = 0,\; P_{\alpha}^{\, T}\, (H_{\overline{\mathcal{E}}} + H_{\overline{\mathcal{E}}}^T)\, P = 0\right\} = \{ 0 \}.
	\end{equation}
	Moreover, a necessary condition for $\overline{X}$ to be constraint nondegenerate is
	$$
	(n-|\alpha|)(n-|\alpha|+1)/2  \leq |\mathcal{E}|.
	$$
\end{proposition}
\proof
	One can easily check that the condition \eqref{LICQ} can be  rewritten as
	$$
	{\rm lin}\left(\mathcal{T}_{\,\mathbb{S}_+^n}\left(\,\overline{X}\,\right)\right)  + {\rm lin}\left(\mathcal{T}_{\,\mathbb{N}^n}\left(\,\overline{X}\,\right)\right) = \mathbb{S}^n,
	$$
	or equivalently,
	$${\rm lin}\left(\mathcal{T}_{\,\mathbb{S}_+^n}\left(\,\overline{X}\,\right)\right) ^\perp \cap {\rm lin}\left(\mathcal{T}_{\,\mathbb{N}^n}\left(\,\overline{X}\,\right)\right)^\perp  = \{0\}.$$
	Direct computation shows that
	$$
	{\rm lin}\left(\mathcal{T}_{\,\mathbb{S}_+^n}\left(\,\overline{X}\,\right)\right)^\perp = \left\{ H \in \mathbb{S}^n \mid P_{\alpha}^{\,T}\, H\,P = 0\right\}, \quad {\rm lin}\left(\mathcal{T}_{\,\mathbb{N}^n}\left(\,\overline{X}\,\right)\right)^\perp=\left\{H\in \mathbb{S}^n\mid H_{{\mathcal{E}}} = 0\right\},
	$$
	which yields the equivalence of \eqref{prop:licq:eqiv} and the definition of  constraint nondegeneracy in \eqref{LICQ}. To complete the proof of this proposition,  we observe that
	$$
	\left\{\begin{array}{ll}
	{\rm dim}\left({\rm lin}\left(\mathcal{T}_{\,\mathbb{S}_+^n}\left(\,\overline{X}\,\right)\right)\right) = n(n+1)/2 - (n-|\alpha|)(n-|\alpha|+1)/2, \\[0.2in]
	{\rm dim}\left({\rm lin}\left(\mathcal{T}_{\,\mathbb{N}^n}\left(\,\overline{X}\,\right)\right)\right) = |\,\mathcal{E}\,|,
	\end{array}\right.
	$$
	where dim$(S)$ represents the dimension of a given linear space $S$. Therefore, a necessary condition for the constraint nondegeneracy to hold at a feasible point $\overline{X}$ is
	$$
	n(n+1)/2 - (n-|\alpha|)(n-|\alpha|+1)/2 + |\,\mathcal{E}\,|\, \geq\,  n(n+1)/2.
	$$
	From here, the required result follows.
\endproof

\begin{remark}\label{remark:LICQ}
	Proposition \ref{prop:LICQ} indicates that the feasible point $\overline{X}$ is likely to be degenerate if either the rank of $\overline{X}$ or the number of nonzero entries of $\overline{X}$ are small.
\end{remark}

In the following, we discuss a consequence of the constraint degeneracy to the Newton-type algorithm for solving \eqref{eq:DNN}. Observe that the Slater condition always holds for \eqref{eq:DNN}, which implies the existence of  optimal solutions to \eqref{dual} \cite[Theorem 28.2]{rockafellar1970convex}. Moreover, the unique optimal solution  $\overline{X}\in \mathbb{S}^n$ to \eqref{eq:DNN} and any dual optimal solution $(\overline{S}, \overline{Z})\in\mathbb{S}^n\times \mathbb{S}^n$ form a KKT point to \eqref{eq:DNN}, at which $\Res(\overline{X}, \overline{S}, \overline{Z}) = 0$ \cite[Theorem 28.3]{rockafellar1970convex}, where $\Res(\cdot)$ is the KKT residual function defined in \eqref{defn:kkt}. Notice that the function $\Res(\cdot)$ is globally Lipschitz continuous so that it is F(r{\'e}chet)-differentiable almost everywhere \cite[Section 9.J]{RRWets1998}. This fact makes the following Clarke's generalized Jacobian of $\Res$ at any $(X,S,Z)\in\mathbb{S}^n\times \mathbb{S}^n\times \mathbb{S}^n$ well defined:
$$
\partial\, \Res(X,S,Z) \,\triangleq \, {\rm conv}\{\partial_B \, \Res(X,S,Z)\},
$$
where for any $W = (X,S,Z)$,
$$
\partial_B \, \Res(W)\,\triangleq \,\left\{V\in \mathbb{S}^n\times \mathbb{S}^n\times \mathbb{S}^n\,\mid\, V = \lim_{k\to \infty} \Res^{\,\prime}(W^k),\; W^k \to W,\; \Res\; {\rm is}\; {\rm F}\text{-}{\rm differentiable \; at}\; W^k \right\}.
$$
Moreover, the function $\Res(\cdot)$ is strongly semismooth since both $\Pi_{\mathbb{S}_+^n}(\cdot)$ \cite{Sun2002semismooth} and $\Pi_{\mathbb{N}^n}(\cdot)$ \cite[Proposition 7.4.7]{facchinei2007finite} are strongly semismooth. Thus the semismooth Newton method can be applied to solve the semismooth equation $\Res(X,S,Z) = 0$, where the $(k+1)$-th Newton direction $d\in\mathbb{S}^n\times \mathbb{S}^n\times \mathbb{S}^n$ is the solution of the following linear equation (c.f. \cite[Section 7.5]{facchinei2007finite}):
$$
\Res(X^{\,k},S^{\,k},Z^{\,k}) + V^k d = 0, \quad V^k\in \partial\, \Res(X^{\,k},S^{\,k},Z^{\,k}).
$$
Though the local superlinear convergence of this method can be established under the nonsingularity of $\partial\, \Res(\overline{X}, \overline{S}, \overline{Z})$ at a KKT point $(\overline{X}, \overline{S}, \overline{Z})$, the following proposition however reveals that such a nonsingularity condition cannot hold if $\overline{X}$ is constraint degenerate.
\begin{proposition}\label{prop:LICQ:Jacobian}
	Let $\overline{X}\in\mathbb{S}^n$ be the unique optimal solution to \eqref{eq:DNN}. Let $(\overline{S}, \overline{Z})\in \mathbb{S}^n\times \mathbb{S}^n$ be an optimal solution to \eqref{dual} such that $(\overline{X}, \overline{S}, \overline{Z})$ is a KKT point of \eqref{eq:DNN}. Then any element in $\partial \,\Res(\overline{X}, \overline{S}, \overline{Z})$ is nonsingular if and only if $\overline{X}$ is constraint nondegenerate.
\end{proposition}
\proof
	It is known from
	\cite[Theorem 4.1]{SunDF2006strong} that for a general nonlinear semidefinite programming problem, which includes \eqref{eq:DNN} as a special case, any element in $\partial \Res(\overline{X}, \overline{S}, \overline{Z})$ is nonsingular if and only if the strong second order sufficient condition holds at $\overline{X}$ and $\overline{X}$ is constraint nondegenerate. Since the objective function in \eqref{eq:DNN} is strongly convex, the strong second order sufficient condition obviously holds at $\overline{X}$. Therefore, the conclusion of this proposition follows.
\endproof

Based on Proposition \ref{prop:LICQ:Jacobian}, we see that it is not suitable to adopt the semismooth Newton method to solve the equation \eqref{defn:kkt} if the optimal solution $\overline{X}$ of \eqref{eq:DNN} is degenerate. According to Remark \ref{remark:LICQ}, this degeneracy is likely to occur when $\overline{X}$ has low rank or is sparse, a situation that may be frequently encountered in practical applications. To avoid such an unfavorable situation for the semismooth Newton method, we  design an ALM in the next section for solving the problem \eqref{eq:DNN}, for which the semismooth Newton method is employed to solve a sequence of nonsingular semismooth equations.

\section{The Dual Quadratic Growth Condition and the Asymptotically Superlinear Convergence of the ALM}\label{sec:growth condition}

In this section, we first take a detour to discuss sufficient conditions for the quadratic growth condition of a general class of convex constrained optimization problems, which includes \eqref{eq:DNN} as a special case. These sufficient conditions will be used to derive the asymptotically superlinear convergence rate of the ALM in \eqref{ALM:iter}, to be presented in the last part of this section.

\subsection{Sufficient conditions for the quadratic growth condition}

Let $F:\X\rightrightarrows\Y$ be a set-valued mapping. The graph of the mapping $F$ is defined as $\text{gph}\,(F) \,\triangleq \, \{\,(x,y)\in\X\times \Y\,\mid\,y\in F(x)\,\}$. The following definition  of metric subregularity  is taken from~\cite[Section 3.8(3H)]{dontchev2009implicit}.

\begin{definition}\label{defn:metric11}
	A set-valued mapping  $F:\X\rightrightarrows\Y$ is said to be metrically subregular at $\bar x\in\X$ for $\bar y\in\Y$ with modulus $\kappa > 0$ if $(\bar x, \bar y)\in \textup{gph}\,(F)$ and there exist a constant $\varepsilon>0$ such that
	\[
	{\rm dist}\left(\,x, F^{-1}(\bar y)\,\right)\, \leq \, \kappa\,{\rm dist}\left(\,\bar y, F(x)\,\right), \quad\forall\, x\in\mathbb{B}_{\varepsilon}(\bar{x}).
	\]
\end{definition}

The next result, which provides a convenient way to check the metric subregularity of the subdifferential of a proper closed convex function, is proven in~\cite[Theorem 3.3]{aragon2008characterization}.

\begin{proposition}\label{pre:metricregular:thm}
	Let $\mathcal{H}$ be a real Hilbert space endowed with the inner product $\langle \cdot,\cdot\rangle$ and $\theta:\mathcal{H}\to(-\infty, +\infty]$ be a proper lower semicontinuous convex function. Consider $\bar{x}, \bar{v}\in\mathcal{H}$ satisfying $(\bar{x},\bar{v})\in\textup{gph}\,(\partial \theta)$. Then $\partial \theta$ is metrically subregular at $\bar{x}$ for $\bar{v}$ if and only if there exist	constants $\kappa > 0$ and $\varepsilon > 0$ such that
	\[
	\theta(x)\, \geq \, \theta(\bar{x}) + \langle \bar{v}, x-\bar{x}\rangle  + \kappa\,{\rm dist}^2\left(\,x, (\partial \theta)^{-1}(\bar{v})\,\right), \quad\forall\, x\in\mathbb{B}_{\varepsilon}(\bar{x}).
	\]
\end{proposition}

A set-valued mapping ${F}:\X\rightrightarrows\Y$ is said to be  polyhedral if its graph is the union of finitely many polyhedral convex sets. Below is a fundamental result from Robinson~\cite{robinson1981some} on polyhedral mappings.

\begin{proposition}\label{prop:polyhedral}
	Let ${F}:\X\rightrightarrows\Y$ be a set-valued polyhedral mapping and $(\bar{x},\bar{y})\in\textup{gph}(F)$. Then $F$ is locally upper Lipschitz continuous at $\bar{x}$, i.e., there exist constants $\kappa>0$ and $\varepsilon>0$ such that
	$$
	F(x) \, \subseteq \, F(\bar x) + \kappa\,\|x-\bar x\|\,\mathbb{B}_{1}(0), \quad\forall\, x\in\mathbb{B}_{\varepsilon}(\bar{x}).
	$$
\end{proposition}

In our subsequent discussions, we also need the concept of bounded linear regularity of a collection of closed convex sets, which can be found from, e.g., ~\cite[Definition 5.6]{bauschke1996projection}.
\begin{definition}
	Let $D_1, D_2, \ldots, D_m\subseteq\X$ be closed convex sets for some positive integer $m$. Suppose that $D \,\triangleq \, D_1\cap D_2 \cap\ldots\cap D_m$ is non-empty. The collection $\{D_1, D_2, \ldots, D_m\}$ is said to be boundedly linearly regular if for every bounded set $B\subseteq \X$, there exists a constant $\kappa >0$ such that
	$$
	{\rm dist}\left(\,x, D\,\right) \, \leq \, \kappa\, \max\left\{{\rm dist}\left(\,x, D_1\,\right), \ldots, {\rm dist}\left(\,x, D_m\,\right)\right\}, \;\forall\, x\in B.
	$$
\end{definition}

A sufficient condition to guarantee the property of bounded linear regularity is established  in~\cite[Corollary 3]{bauschke1999strong}.

\begin{proposition}\label{prop:boundedlinear}
	Let $D_1, D_2, \ldots, D_m\subseteq\X$ be closed convex sets for some positive integer $m$. Suppose that $D_1, D_2, \ldots, D_r$ are polyhedral for some $r\in\{0,1,\ldots, m\}$. Then a sufficient condition for $\{D_1, D_2, \ldots, D_m\}$ to be boundedly linearly regular is
	$$
	\bigcap_{i=1,2, \ldots, r} D_i \quad \cap\;  \bigcap_{i = r+1, \ldots, m} \textup{ri}\,(D_i)\neq \emptyset.
	$$
\end{proposition}

Consider the following linear equality and inequality constrained nonsmooth convex problem:
\begin{equation}\label{eb:opt}
\begin{array}{cl}
\displaystyle\operatornamewithlimits{minimize}_{x\in \mathbb{X}} & \theta(x) \,\triangleq \,  h(\mathcal{F}x) +\langle c,x\rangle + p(x) \\[3pt]
\mbox{subject to} & \mathcal{A}x -b \in \mathcal{Q},
\end{array}
\end{equation}
where $\mathcal{F}:\X\to\W$ and $\mathcal{A}:\X\to\Y$ are given linear maps, $\mathcal{Q}\subseteq\mathbb{Y}$ is a given convex polyhedral cone,	$c\in \X$ and $b\in\Y$ are given data, $p:\X\to (-\infty, +\infty]$ is a closed proper convex  function, $h:\W\to (\infty, +\infty]$ is an essentially smooth and essentially locally strictly convex function. The Lagrangian dual of this problem is
\begin{equation}
\begin{array}{ll}
\displaystyle\operatornamewithlimits{maximize}_{y\in \mathbb{Y}}\, & g(y)\,\triangleq \, \displaystyle\inf_{x\in\X}\, \{\, \theta(x) + \langle y, \mathcal{A}x - b\rangle\, \}, \\[3pt]
\mbox{subject to} & y\in\mathcal{Q}^\circ.
\end{array}
\label{eb:dual}
\end{equation}
Assume that the following KKT system associated with problem \eqref{eb:opt} admits at least one solution:
\begin{equation}\label{kkt}
\left\{\begin{array}{ll}
0 \in \mathcal{F}^*\nabla h(\mathcal{F}x)  +c +  \partial p(x) + \mathcal{A}^*y,\\[5pt]
y\in \mathcal{N}_{\mathcal{Q}}(\mathcal{A}{x} - b),
\end{array}\right. \quad (x,y)\in\X\times \Y.
\end{equation}
We denote ${\rm SOL}_{\rm P}$ as the solution set of problem \eqref{eb:opt} and ${\rm SOL}_{\rm D}$ as the solution set of problem \eqref{eb:dual}. It is known from~\cite[Theorem 30.4 \text{and} Corollary 30.5.1]{rockafellar1970convex} that $(\bar{x}, \bar{y})\in\X\times \Y$ solves the KKT system (\ref{kkt}) if and only if $\bar{x}\in {\rm SOL}_{\rm P}$  and $\bar{y}\in {\rm SOL}_{\rm D}$. To further characterize  ${\rm SOL}_{\rm P}$, we need the following invariant property of $\mathcal{F}x$ over  ${\rm SOL}_{\rm P}$, whose proof readily follows from  well-known existing techniques in the literature \cite{luo1992linear,mangasarian1988simple,tseng2010approximation}.

\begin{lemma}\label{lemma:invariant}
	The value $\mathcal{F}{x}$ is invariant over $x\in {\rm SOL}_{\rm P}$, i.e.,
	for any $x', x''\in{\rm SOL}_{\rm P}$, we have $\mathcal{F}{x'} = \mathcal{F}x''$.
\end{lemma}

Take an arbitrary point $\bar{x}\in{\rm SOL}_{\rm P}$ and denote
\begin{equation}\label{defn:ubar}
\begin{array}{ll}
\bar{\zeta}\,\triangleq\, \mathcal{F}\bar{x}, \quad
\bar{\eta} \,\triangleq\, \mathcal{F}^*\nabla h(\bar{\zeta}) +c,\quad
\overline{\mathcal{V}}\,\triangleq \,\{x\in\X\; {\mid}\;  \mathcal{F}x = \bar{\zeta}\}.  \end{array}
\end{equation}
We define two set-valued mappings $\mathcal{G}_1:\Y\rightrightarrows\X$ and $\mathcal{G}_2:\Y\rightrightarrows\X$ by
\begin{equation}\label{defn:mappingG}
\mathcal{G}_1(y) \,\triangleq \, (\partial p)^{-1}(-\mathcal{A}^*y - \bar{\eta}),\quad  \quad
\mathcal{G}_2(y) \,\triangleq \, \{x\in\mathbb{X}\, {\mid}\, y\in \mathcal{N}_{\mathcal{Q}}(\mathcal{A}x -b)\}, \quad y\in\Y.
\end{equation}
Then, from (\ref{kkt}), Lemma \ref{lemma:invariant} and the discussion above Lemma \ref{lemma:invariant}, we immediately obtain the following useful observation for the optimal solution set ${\rm SOL}_{\rm P}$.

\begin{proposition}\label{prop:Omega}
	Assume that $\bar{x}\in {\rm SOL}_{\rm P}$ and $\bar{y}\in {\rm SOL}_{\rm D}$. Then the  optimal solution set ${\rm SOL}_{\rm P}$ can be characterized as
	$$\begin{array}{ll}
	{\rm SOL}_{\rm P}  = \{x\in\X \; {\mid} \; \mathcal{F}x = \bar{\zeta}, \; 0\in \bar{\eta} + \partial p({x})
	+ \mathcal{A}^*\bar{y},
	\; \bar{y}\in \mathcal{N}_{\mathcal{Q}}(\mathcal{A}x -b)\} = \overline{\mathcal{V}}\cap \mathcal{G}_1(\bar{y})\cap \mathcal{G}_2(\bar{y}).
	\end{array}
	$$
\end{proposition}

The following concept of quadratic growth condition for problem (\ref{eb:opt}) plays an important role in our later analysis.
\begin{definition}
	The quadratic growth condition for problem (\ref{eb:opt}) holds at an optimal solution $\bar{x}\in {\rm SOL}_{\rm P}$ if there exist positive constants $\kappa$ and $\varepsilon$ such that
	\begin{equation}\label{defn:quadraticGrowth}
		\theta(x) \geq \theta(\bar{x}) + \kappa\,{\rm dist}^2(x,  {\rm SOL}_{\rm P}), \quad\forall\; x\in\mathbb{B}_{\varepsilon}(\bar{x})\cap\{x\in\X\,\mid\,  \mathcal{A}x-b \in\mathcal{Q}\,\}.
	\end{equation}
\end{definition}

To analyze the quadratic growth condition of problem \eqref{eb:opt}, we will need the following assumption and lemma later.

\begin{assumption}\label{ass:f}
	The following local growth conditions hold:\\[3pt]
	(a) For any ${w}\in\textup{dom}\;h$, there exist positive constants $\kappa_1$  and $\varepsilon_1$ such that
	$$
	h(w')\geq h({w}) + \langle \nabla h({w}), w' - {w}\rangle  + \kappa_1\|w'-{w}\|^2, \quad\forall\, w'\in\mathbb{B}_{\varepsilon_1}(w).
	$$
	(b) For any $(x,v)\in\textup{gph}\, (\partial p)$, there exist positive constants $\kappa_2$  and $\varepsilon_2$ such that
	$$
	p(x')\geq p(x) + \langle v, x' - {x}\rangle  + \kappa_2\,{\rm dist}^2\left(x', (\partial p)^{-1}(v)\right), \quad\forall \, x'\in\mathbb{B}_{\varepsilon_2}(x).
	$$
\end{assumption}

\begin{lemma}\label{lemma:polyhdist}
	Let $\bar{x}\in{\rm SOL}_{\rm P}$ and $\bar{y}\in{\rm SOL}_{\rm D}$. Then there exist positive constants $\kappa$ and $\varepsilon$ such that
	$$
	{\rm dist}\,(x, \mathcal{G}_2(\bar{y}))\leq \kappa \, {\rm dist}\,(\mathcal{A}x-b, \mathcal{N}_{\mathcal{Q}^\circ}(\bar{y})), \quad\forall\, x\in\mathbb{B}_{\varepsilon}(\bar{x}).
	$$
\end{lemma}
\proof
First we note that since $\mathcal{Q}$ is a closed convex cone, $y\in \mathcal{N}_{\mathcal{Q}}(z)$ if and only if $z\in \mathcal{N}_{\mathcal{Q}^\circ}(y)$. Define the subspace $\Xi_1\subseteq\X\times \Y$ and {the} polyhedral set $\Xi_2\subseteq\X\times \Y$ by 
$$ 
\Xi_1 = \{(x,q)\in\X\times \Y\,\mid\, \mathcal{A}x-b = q\},\quad\quad \Xi_2 = \{(x,q)\in\X\times \Y\,\mid\, q\in\mathcal{N}_{\mathcal{Q}^\circ}(\bar{y})\}.
$$
Denote $\widetilde{\mathcal{G}}_2\,\triangleq \,\Xi_1\cap\Xi_2$, which is non-empty as $(\bar{x},\mathcal{A}\bar{x}-b)\in\widetilde{\mathcal{G}}_2$. Since $\Xi_1$ and $\Xi_2$ are polyhedral sets, we know from Proposition \ref{prop:boundedlinear} that the collection $\{\Xi_1, \Xi_2\}$ is boundedly linearly regular.  Therefore, there exist positive constants $\kappa$ and $\varepsilon$ such that for any $x\in \mathbb{B}_{\varepsilon}(\bar{x})$,
$$
\text{dist}\,((x,\mathcal{A}x-b), \widetilde{\mathcal{G}}_2)\leq
\kappa\big(\text{dist}\,((x,\mathcal{A}x-b),\Xi_1) +
\text{dist}\,((x,\mathcal{A}x-b),\Xi_2)\big) =  \kappa\,\text{dist}\,(\mathcal{A}x-b,\mathcal{N}_{\mathcal{Q}^\circ}(\bar{y})).
$$
Now note that there exists $(x',w')\in\widetilde{\mathcal{G}}_2$ such that
$$\text{dist}((x,\mathcal{A}x-b), \widetilde{\mathcal{G}}_2) = \sqrt{\|x - x'\|^2 + \|\mathcal{A}x -b- w'\|^2}\geq \|x-x'\|\geq
\text{dist}(x, \mathcal{G}_2(\bar{y})),
$$
where the last inequality follows from the fact that
$x'\in \mathcal{G}_2(\bar{y})$ because $\mathcal{A}x'-b\in \mathcal{N}_{\mathcal{Q}^\circ}(\bar{y})$ implies that $\bar{y} \in\mathcal{N}_{\mathcal{Q}}(\mathcal{A} x' -b).$  From here, we complete the proof of the lemma.
\endproof

The following result, which is partially motivated by the recent paper~\cite{zhou2015eb} and its  further development in~\cite{drusvyatskiy2016error} for convex composite  optimization problems regularized by the nuclear norm function of rectangular matrices, provides a general approach for proving the quadratic growth condition of (\ref{eb:opt}) where the constraint $\mathcal{A}x-b\in\mathcal{Q}$ is present.

\begin{theorem}\label{prop:quadraticgrow}
	Assume that ${\rm SOL}_{\rm P}$  is non-empty. Suppose that Assumption \ref{ass:f} holds and that there exists $\bar{y}\in{\rm SOL}_{\rm D}$ such that the collection of three sets $\{\overline{\mathcal{V}}, \mathcal{G}_1(\bar{y}),\mathcal{G}_2(\bar{y})\}$ is boundedly linearly regular.	Then the  quadratic growth condition (\ref{defn:quadraticGrowth}) holds at any  $\bar{x}\in{\rm SOL}_{\rm P}$.
\end{theorem}

\proof
Let $\bar{x}\in{\rm SOL}_{\rm P}$ be an arbitrary but fixed point. Since $(\bar{x},-\mathcal{A}^* \bar{y}-\bar{\eta}) \in {\rm gph}(\partial p)$, from Assumption \ref{ass:f} (b), we know that there exist positive constants $\kappa_1$ and  $\varepsilon$ such that
$$
p(x) \, \geq\,  p(\bar{x}) + \left\langle\, -\mathcal{A}^*\bar{y} - \bar{\eta}, x - \bar{x}\,\right\rangle  + \kappa_1 \, {\rm dist}^2\left(x, (\partial p)^{-1}(-\mathcal{A}^*\bar{y} - \bar{\eta})\right), \;\forall \,x\in\mathbb{B}_{\varepsilon}(\bar{x}).
$$
Note that $(\mathcal{A}\bar{x}-b,\bar{y})\in {\rm gph}\left(\mathcal{N}_{\mathcal{Q}^\circ}^{-1}\right)$ and  $\mathcal{N}_{\mathcal{Q}^\circ}(\cdot)$ is a set-valued polyhedral function. Also, $\mathcal{N}_{\mathcal{Q}^\circ}^{-1} = \partial \delta_{\mathcal{Q}}.$ Thus, we can obtain from Proposition \ref{prop:polyhedral} that $\mathcal{N}_{\mathcal{Q}^\circ}(\cdot)$ is locally upper Lipschitz continuous, which further implies the metric subregularity of $\mathcal{N}_{\mathcal{Q}^\circ}^{-1}$  at $\mathcal{A}\bar{x}-b$ for  $\bar{y}$ by definition. Now by shrinking $\varepsilon$ if necessary, we know that there exists a constant {$\kappa_1'>0$} such that
$$
\delta_{\mathcal{Q}}(\mathcal{A}x-b)\geq \delta_{\mathcal{Q}}(\mathcal{A}\bar{x}-b) + \langle\, \bar{y}, \mathcal{A}x -b- (\mathcal{A}\bar{x}-b)\,\rangle  + \kappa'_1\, {\rm dist}^2\left(\mathcal{A}x-b, \mathcal{N}_{\mathcal{Q}^\circ}(\bar{y})\right), \;\forall \,x\in\mathbb{B}_{\varepsilon}(\bar{x}).
$$
Moreover, the assumed bounded linear regularity of $\{\overline{\mathcal{V}}, \mathcal{G}_1(\bar{y}),\mathcal{G}_2(\bar{y})\}$ and  the result in Proposition \ref{prop:Omega} imply that there exist $\kappa_2>0$ and $\kappa_3 >0$, such that for any $x\in\mathbb{B}_{\varepsilon}(\bar{x})$,
$$\begin{array}{ll}
{\rm dist}^2(x,  {\rm SOL}_{\rm P}) &
= {\rm dist}^2\left(x, \overline{\mathcal{V}}\cap\mathcal{G}_1(\bar{y})\cap\mathcal{G}_2(\bar{y})\right)\\[0.15in]
&  \leq
\kappa_2\big[\,{\rm dist}^2(x, \overline{\mathcal{V}}) + {\rm dist}^2(x, \mathcal{G}_1(\bar{y}))+ {\rm dist}^2(x, \mathcal{G}_2(\bar{y}))\,\big]\\[0.15in]
& \leq  \kappa_3\big[\,\|\,\mathcal{F}x - \bar{\zeta}\,\|^2   +
{\rm dist}^2(x, (\partial p)^{-1}(-\mathcal{A}^*\bar{y} -\bar{\eta}))+{\rm dist}^2(\mathcal{A}x-b, \mathcal{N}_{\mathcal{Q}^\circ}(\bar{y}))\,\big],  \end{array}
$$
where in the last inequality, the first term  comes from Hoffman's error bound~\cite{hoffman1952approximate} and the third term comes from Lemma \ref{lemma:polyhdist}. Then by Assumption \ref{ass:f} {(a)}, shrinking $\varepsilon$ if necessary, we know that there exists $\kappa_4>0$ such that for any $x\in\mathbb{B}_{\varepsilon}(\bar{x})$,
$$
\label{ineq:strongconvex}
\begin{array}{ll}
h(\mathcal{F}x) \geq  h(\bar{\zeta}) + \langle\, \nabla h(\bar{\zeta}), \mathcal{F}x - \bar{\zeta}\,\rangle  + \kappa_4\,\|\,\mathcal{F}x - \bar{\zeta}\,\|^2.
\end{array}
$$
Taking all the above inequalities into account and  recalling that $\bar{\eta} = \mathcal{F}^*\nabla h(\bar{\zeta}) + c$ in (\ref{defn:ubar}), we  derive, for any $x\in \mathbb{B}_{\varepsilon}(\bar{x})\cap \{x\in\X\,\mid\,\mathcal{A}x-b\in\mathcal{Q}\}$, that
$$
\begin{array}{rl}
&\theta(x) = h(\mathcal{F} x) + \langle c, x \rangle + p(x) + \delta_{\mathcal{Q}}(\mathcal{A}x-b)
\\[0.15in]
\geq
& \theta(\bar{x}) + \langle\, \mathcal{F}^*\nabla h(\bar{\zeta}) +c- \bar{\eta}, x - \bar{x}\,\rangle +\kappa_4\,\|\,\mathcal{F}x - \bar{\zeta}\,\|^2+ \kappa_1{\rm dist}^2(x, (\partial p)^{-1}(-\mathcal{A}^*\bar{y} - \bar{\eta}))\\[0.15in]
&+ \kappa'_1\, {\rm dist}^2(\mathcal{A}x-b, \mathcal{N}_{\mathcal{Q}^\circ}(\bar{y}))\\[0.15in]
\geq &  \theta(\bar{x})  + \min \{\kappa_1,\kappa_1', \kappa_4\}\left[\,\|\mathcal{F}x - \bar{\zeta}\|^2 + {\rm dist}^2(x, (\partial p)^{-1}(-\mathcal{A}^*\bar{y} - \bar{\eta}))  + {\rm dist}^2(\mathcal{A}x-b, \mathcal{N}_{\mathcal{Q}^\circ}(\bar{y}))\,\right]\\[0.15in]
= &  \theta(\bar{x}) + \kappa_3^{-1}\min \{\kappa_1, \kappa_1',\kappa_4\}\,{\rm dist}^2(x, {\rm SOL}_{\rm P}),
\end{array}
$$
which establishes the desired result.
\endproof

\subsection{The quadratic growth condition for \eqref{dual}}

Notice that \eqref{dual} can be viewed as a special case of \eqref{eb:opt} where $\mathcal{A}$, $b$, $c$ and $\mathcal{Q}$ are vacant, and
\[
\left\{\begin{array}{ll}
h(U) =  \displaystyle\frac{1}{2}\,\left\|\,U+G\,\right\|^2, \;\;  U \in \mathbb{S}^n,\\[0.1in]
\mathcal{F}(S,Z) = S+ Z,
\quad p(S,Z) = \delta_{\mathbb{S}^n_+}(S)
+ \delta_{\mathbb{N}^n}(Z), \;\;  (S,Z)\in \mathbb{S}^n\times \mathbb{S}^n.
\end{array}\right.
\]
In this section, we show that  Assumption \ref{ass:f} always holds for such a case, while the bounded linear regularity of the corresponding sets $\{\overline{\mathcal{V}}, \mathcal{G}_1(\bar{y}),\mathcal{G}_2(\bar{y})\}$ is implied by the existence of a strict complementarity solution of \eqref{eq:DNN}.

Let $\overline{X}\in\S_+^n$ and $\overline{S}\in\S_+^n$ satisfy $0\in\overline{X}+\partial \delta_{\S_+^n}(\overline{S})$, or equivalently, $\langle \overline{X},\overline{S}\rangle  = 0$. Suppose that $\overline{Z}\,\triangleq \,\overline{X} - \overline{S}$ has its eigenvalues $\bar{\lambda}_1\geq \bar{\lambda}_2\geq \ldots\geq \bar{\lambda}_n$ being arranged in {a} non-increasing order. Denote
\begin{equation}\label{defn:indices}
\alpha\,\triangleq \,\{i \mid \bar{\lambda}_i>0, \; 1\leq i\leq n\},\quad
\beta\,\triangleq \,\{i \mid \bar{\lambda}_i=0, \; 1\leq i\leq n\},\quad
\gamma\,\triangleq \,\{i \mid \bar{\lambda}_i<0, \; 1\leq i\leq n\}.
\end{equation}
Then there exists an orthogonal matrix $\overline{P}\in\mathcal{O}^n$ such that
\begin{equation}\label{eig-decom}
\overline{Z} = \overline{P}\left(\begin{array}{ccc}
\overline{\Lambda}_\alpha & & \\
& 0 &\\
& & -\overline{\Lambda}_{\gamma}
\end{array}\right)\overline{P}^{\,T}, \quad
\overline{X}= \overline{P}\left(\begin{array}{ccc}
\overline{\Lambda}_\alpha & & \\
& 0 &\\
& & 0_{|\gamma|}
\end{array}\right)\overline{P}^{\,T},
\quad  \overline{S}= \overline{P}\left(\begin{array}{ccc}
0_{|\alpha|} & & \\
& 0 &\\
& & \overline{\Lambda}_\gamma
\end{array}\right)\overline{P}^{\,T},
\end{equation}
where $\overline{\Lambda}_\alpha = {\rm diag}(\bar{\lambda}_\alpha) \succ 0$ and $\overline{\Lambda}_\gamma = {\rm diag}(-\bar{\lambda}_\gamma) \succ 0$. Denote $\overline{P} = [\,\overline{P}_{\alpha}\; \overline{P}_{\beta}\; \overline{P}_\gamma\,]$ with $\overline{P}_\alpha\in\R^{n\times |\alpha|}$, $\overline{P}_\beta\in\R^{n\times |\beta|}$ and $\overline{P}_\gamma\in\R^{n\times |\gamma|}$. Then we have
$$
\left\{\begin{array}{rll}
\mathcal{T}_{\S_+^n}(\overline{X}) &= &\left\{\,  H\in\S^n \,\mid\, [\,\overline{P}_{\beta}\;\overline{P}_{\gamma}\,]^{\,T}H\,[\,\overline{P}_{\beta}\;\overline{P}_{\gamma}\,]\succeq 0\,\right\}, \\[0.15in]
\mathcal{N}_{\S_+^n}(\overline{X})  &= &\left\{\, H\in\S^n \,\mid\, [\,\overline{P}_{\beta}\;\overline{P}_{\gamma}\,]^{\,T}H\,[\,\overline{P}_{\beta}\;\overline{P}_{\gamma}\,]\preceq 0, \;\overline{P}_{\alpha}^{\,T}H\overline{P}=0\,\right\}.
\end{array}\right.
$$
By noting that $\partial \delta_{\S_+^n}(\overline{S}) = \mathcal{N}_{\S_+^n}(\overline{S})$, we immediately obtain the following results.

\begin{proposition}\label{prop:boundedlinearsdp}
	Let $\overline{S}\in\S_+^n$ and $0\in\overline{X}+\partial \delta_{\S_+^n}(\overline{S})$. Suppose that $\overline{S}$ and $\overline{X}$ have eigenvalue decompositions as in (\ref{eig-decom}). Then it holds that:  \\[5pt]
	(a) $\mathcal{N}_{\S_+^n}(\overline{S})$ is a polyhedral set if and only if $|\gamma| \geq n-1$; \\[5pt]
	(b) $0\in \overline{X} + {\rm ri}\,\left(\mathcal{N}_{\S_+^n}(\overline{S})\right)$ if and only if $|\beta| = 0$, i.e., $\textup{rank}(\overline{X}) + \textup{rank}(\overline{S}) = n$.
\end{proposition}

The following proposition shows that $\mathcal{N}_{\mathbb{S}_+^n}(\cdot)$ is metrically subregular at any point on its graph. This result is part of the first author's PhD thesis~\cite[Section 2.5.2]{CuiThesis}, which can also be derived from the recent work \cite{CDZhao2016}.  However, here we furnish a direct proof for better understandings of the nonpolyhedral semidefinite cone.

\begin{proposition}\label{sdp:metricsub}
	Let $\overline{S}\in\S_+^n$ and $0\in \overline{X} +  \partial \delta_{\S_+^n}(\overline{S})$. Then $\partial \delta_{\S_+^n}(\cdot)$ is metrically subregular at $\overline{X}$ for $-\overline{S}$ and $\partial \delta_{\S_-^n}(\cdot)$ is metrically subregular at $-\overline{S}$ for $\overline{X}$.
\end{proposition}
\proof
In the following, 	we shall prove the metric subregularity of  $\partial \delta_{\S_-^n}(\cdot)$  at $-\overline{S}$ for $\overline{X}$ and its counterpart regarding $\partial\delta_{\mathbb{S}_+^n}$ can be obtained similarly. Without loss of generality, let  $\overline{X}$ and $\overline{S}$ have the eigenvalue decompositions as in (\ref{eig-decom}). According to Proposition \ref{pre:metricregular:thm}, in order to prove the {metric subregularity} of $\partial \delta_{\S_-^n}(\cdot)$  at $-\overline{S}$ for $\overline{X}$, it suffices to show that there exist a constant $\kappa>0$ and a neighborhood $\mathcal{U}$ of $\overline{S}$ such that for any $S\;\in\S_+^n\cap \mathcal{U}$,
\begin{equation}\label{pre:im:metric}
0\geq \langle\, \overline{X}, -S +\overline{S}\,\rangle + \kappa\,{\rm dist}^2\left(-S, \,(\partial\delta_{\S_-^n})^{-1}(\overline{X})\right)=\langle\, \overline{X}, -S +\overline{S}\,\rangle + \kappa\,{\rm dist}^2\left(-S,\, \mathcal{N}_{\S_+^n}(\overline{X})\right).
\end{equation}
If $|\alpha| = 0$, then $\overline{X}=0$ and the inequality (\ref{pre:im:metric}) holds automatically for any $\kappa\geq 0$ and any neighborhood $\mathcal{U}$ of $\overline{S}$. Thus, we only need to consider the case that $|\alpha|\not=0$. Since the case that $|\gamma|=0$ can be proved similarly as in the case for $|\gamma|\neq 0$, {we} only consider the latter case. Set $\rho \,\triangleq \, \min\{|\bar{\lambda}_j| \mid j\in \gamma\}>0$. Let $S\in \S_+^n\cap\mathbb{B}_{\rho}(\overline{S})$ be arbitrarily chosen. We write $\widetilde{S} = \overline{P}^{\,T}S\overline{P}$ and decompose $\widetilde{S}$ into the following form:
$$
\widetilde{S}\equiv \left(\begin{array}{ccc}
\widetilde{S}_{\alpha\alpha} &\widetilde{S}_{\alpha\beta}  & \widetilde{S}_{\alpha\gamma}\\[4pt]
\widetilde{S}_{\alpha\beta}^{\,T} &\widetilde{S}_{\beta\beta}  & \widetilde{S}_{\beta\gamma}\\[4pt]
\widetilde{S}_{\alpha\gamma}^{\,T} &\widetilde{S}_{\beta\gamma}^{\,T}  & \widetilde{S}_{\gamma\gamma}
\end{array}\right).
$$
By the fact that $S\in\S_+^n$, we can easily check that
$$
\Pi_{\mathcal{N}_{\S_+^n}(\overline{X})}(-S)= -\overline{P}\left(\begin{array}{ccc}
0 & 0 & 0\\[3pt]
0 &\widetilde{S}_{\beta\beta}  & \widetilde{S}_{\beta\gamma} \\[3pt]
0 & \widetilde{S}_{\beta\gamma}^{\, T}   & \widetilde{S}_{\gamma\gamma} \end{array}\right)\overline{P}^{\,T}.
$$
Thus
	\begin{eqnarray}
	\inprod{\overline{X}}{-S + \overline{S}} = \inprod{\overline{X}}{-S} = \inprod{\overline{\Lambda}_\alpha}{-\widetilde{S}_{\alpha\alpha}}
	\leq -\bar{\lambda}_{|\alpha|} {\rm tr}(\widetilde{S}_{\alpha\alpha}).
	\label{eq-temp}
	\end{eqnarray}
In addition, we have
\begin{equation}\label{pre:metricregular:ineq00}
{\rm dist}^2\left(-S, \mathcal{N}_{\S_+^n}(\overline{X})\right)
=\left\|-S - \Pi_{\mathcal{N}_{\S_+^n}(\overline{X})}(-S)\right\|^2
=   \left\|\,\widetilde{S}_{\alpha\alpha}\,\right\|^2+
2\left\|\,\widetilde{S}_{\alpha\beta}\,\right\|^2 + 2\left\|\,\widetilde{S}_{\alpha\gamma}\,\right\|^2.
\end{equation}
Next we proceed to estimate $\norm{\widetilde{S}_{\alpha\alpha}},\norm{\widetilde{S}_{\alpha\beta}}$ and $\norm{\widetilde{S}_{\alpha\gamma}}$. By using the Bauer-Fike Theorem \cite{BauerFike60}, one obtains that for any $i=1,\dots,|\gamma|$,
$$
\begin{array}{ll}
{\rm dist}\left(\lambda_{i}(\widetilde{S}_{\gamma\gamma}), \{\,|\bar{\lambda}_j| \mid j\in \gamma\}\right)
\, \leq \, \left\|\,\widetilde{S}_{\gamma\gamma} -  \overline{\Lambda}_\gamma\,\right\|
\, =\, \left\|\,\overline{P}_\gamma^{\,T} S \overline{P}_\gamma - \overline{P}_\gamma^{\,T}\, \overline{S}\, \overline{P}_\gamma\,\right\|
\, \leq \, \left\|\,S - \overline{S} \,\right\| \, \leq \, \rho.
\end{array}
$$
The above inequality further implies that $0 < \lambda_{i}(\widetilde{S}_{\gamma\gamma})\leq
	|\bar{\lambda}_{n}| + \rho$ for all $i=1,\ldots,|\gamma|$. Thus, $\widetilde{S}_{\gamma\gamma}$ is positive definite and $ \lambda_{\max}(\widetilde{S}_{\gamma\gamma}) \leq  |\bar{\lambda}_{n}|+ \rho.$ Note that $\norm{\widetilde{S}_{\alpha\alpha}} \leq \rho$ and $\norm{\widetilde{S}_{\beta\beta}} \leq \rho$ since $S\in\mathbb{B}_{\rho}(\overline{S})$. Moreover, $\norm{\widetilde{S}_{\alpha\alpha}}^2 \leq \rho \norm{\widetilde{S}_{\alpha\alpha}} \leq \rho\, {\rm tr}(\widetilde{S}_{\alpha\alpha})$.

Now, from the fact that $\widetilde{S}_{\alpha\alpha} - \widetilde{S}_{\alpha\gamma}\,\widetilde{S}_{\gamma\gamma}^{\, -1}\,\widetilde{S}^{\,T}_{\alpha\gamma} \,\succeq \,0$
(because $\widetilde{S} \in\mathbb{S}_+^n$), we have
\begin{equation*}
\lambda_{\textup{max}}^{-1} (\widetilde{S}_{\gamma\gamma})\, \widetilde{S}_{\alpha\gamma}\,\widetilde{S}_{\alpha\gamma}^{\,T}
\;\preceq\;
\widetilde{S}_{\alpha\gamma}\,\widetilde{S}_{\gamma\gamma}^{-1}\,\widetilde{S}_{\alpha\gamma}^{\,T}
\;\preceq\; \widetilde{S}_{\alpha\alpha} .
\end{equation*}
Hence,
\begin{equation} \label{pre:metricregular:ineq11}
\left\|\,\widetilde{S}_{\alpha\gamma}\,\right\|^2 \, = \, \text{tr}\left(\widetilde{S}_{\alpha\gamma}\,\widetilde{S}_{\alpha\gamma}^{\,T}\right)
\, \leq \, \text{tr}\left(\widetilde{S}_{\alpha\alpha}\right)\;\lambda_{\textup{max}}\left(\widetilde{S}_{\gamma\gamma}\right)
\;\leq\;
(|\bar{\lambda}_{n}| + \rho)\, {\rm tr}(\widetilde{S}_{\alpha\alpha}).
\end{equation}
Moreover, we obtain from
$
\left(\begin{array}{cc}
\widetilde{S}_{\alpha\alpha} & \widetilde{S}_{\alpha\beta} \\[0.1in]
\widetilde{S}_{\alpha\beta}^{\,T} & \widetilde{S}_{\beta\beta}
\end{array}\right)\succeq 0
$ that
$$
\widetilde{S}_{ij}^{\,2}\, \leq \, \widetilde{S}_{ii}\,\widetilde{S}_{jj}\, \leq \, \rho\,\widetilde{S}_{ii}
\quad \forall\; i\in\alpha,\, j\in\beta,
$$
which implies that
\begin{equation}\label{pre:metricregular:ineq22}
\left\|\,\widetilde{S}_{\alpha\beta}\,\right\|^2 \, = \, \sum_{i\in\alpha, j\in\beta}\widetilde{S}_{ij}^{\,2} \, \leq \,
\rho \sum_{i\in\alpha, j\in\beta}\widetilde{S}_{ii} \,=\, \rho |\beta| {\rm tr}(\widetilde{S}_{\alpha\alpha}).
\end{equation}
By using the above estimates of $\norm{ \widetilde{S}_{\alpha\alpha}}$, $\norm{\widetilde{S}_{\alpha\beta}}$, and $\norm{\widetilde{S}_{\alpha\gamma}}$ in \eqref{pre:metricregular:ineq00}, we get
	$$
	{\rm dist}^2\left(-S, \mathcal{N}_{\S_+^n}(\overline{X})\right) \;\leq\;
	(2|\bar{\lambda}_n| + 3\rho + 2|\beta|\rho)\, {\rm tr}(\widetilde{S}_{\alpha\alpha}).
	$$
Let $\kappa\,\triangleq \, \displaystyle\frac{\bar{\lambda}_{|\alpha|}}{2|\bar{\lambda}_n| + 3\rho + 2|\beta|\rho}>0$. Then, together with \eqref{eq-temp}, we obtain that $S\in \S_+^n\cap\mathbb{B}_{\rho}\left(\,\overline{S}\,\right)$,
	$$
	\left\langle\, \overline{X}, -S+\overline{S}\,\right\rangle + \kappa\,{\rm dist}^2\left(-S, \, \mathcal{N}_{\S_+^n}(\overline{X})\right)
	\;\leq\; -\bar{\lambda}_{|\alpha|} {\rm tr}(\widetilde{S}_{\alpha\alpha}) + \bar{\lambda}_{|\alpha|} {\rm tr}(\widetilde{S}_{\alpha\alpha})
	=0.
	$$
Therefore, the inequality (\ref{pre:im:metric}) holds for any $S\in \S_+^n\cap\mathbb{B}_{\rho}\left(\overline{S}\right)$ and the proof is completed.
\endproof

Notice that ${\rm SOL}_{\rm P}$ of \eqref{eq:DNN} is a singleton. Combining Theorem \ref{prop:quadraticgrow} and Propositions \ref{prop:boundedlinearsdp} and \ref{sdp:metricsub}, we obtain the following result.
\begin{corollary}\label{coro:sdpTphi}
	Let $\overline{X}\in {\rm SOL}_{\rm P}$ be the unique optimal solution of \eqref{eq:DNN}.
	The quadratic growth condition of \eqref{dual} holds at any $(\overline{S}, \overline{Z})\in {\rm  SOL}_{\rm D}$  under one of the following two conditions:\\[5pt]
	(i)  $\textup{rank}(\overline{X})\geq n-1$;\\[5pt]
	(ii) there exists  $(\widehat{S}, \widehat{Z})\in{\rm SOL}_{\rm D}$ such that ${\rm rank}(\overline{X}) + {\rm rank}(\widehat{S}) = n$.
\end{corollary}
\proof
	Obviously, the function $h(\cdot)= \frac{1}{2}\left\|\,\cdot +G\,\right\|^2$ satisfies Assumption \ref{ass:f}(a) for the point $(\overline{S},\overline{Z})$. For the function $p(S,Z) = \delta_{\mathbb{S}^n_+}(S)
	+ \delta_{\mathbb{N}^n}(Z)$, since $\partial\, p(S,Z) = \mathcal{N}_{\mathbb{S}_+^n}(S) \times \mathcal{N}_{\mathbb{N}^n}(Z)$ and $\mathcal{N}_{\mathbb{N}^n}(\cdot)$ is a polyhedral mapping,  we know from Propositions \ref{pre:metricregular:thm} and \ref{prop:polyhedral} that for any $\overline{V}\in \mathcal{N}_{\mathbb{N}^n}(\overline{Z})$, there exist positive scalars $\varepsilon$ and $\kappa$ such that
	\[
	\delta_{\mathbb{N}^n}(Z) \geq \delta_{\mathbb{N}^n}(\overline{Z}) +  \left\langle \, \overline{V}, Z - \overline{Z}\,\right\rangle + \kappa \, {\rm dist}^2\left(Z,\, \mathcal{N}^{-1}_{\mathbb{N}^n}(\overline{V})\right), \quad \forall \; Z\in \mathbb{B}_\varepsilon(\overline{Z}).
	\]
	This, together with Proposition \ref{sdp:metricsub}, implies Assumption \ref{ass:f}(b) at $(\overline{S},\overline{Z})$. In addition, it is known from Propositions \ref{prop:boundedlinear} and \ref{prop:boundedlinearsdp} that the bounded linear regularity of the polyhedral set $\overline{\mathcal{V}} = \{(S,Z)\in \mathbb{S}^n\times \mathbb{S}^n\mid S+Z = \overline{S} + \overline{Z} \}$ and the nonpolyhedral set	$\mathcal{G}_1(\overline{X}) =  \mathcal{N}^{-1}_{\mathbb{S}_+^n}(\overline{X}) \times \mathcal{N}^{-1}_{\mathbb{N}^n}(\overline{X})$ can be implied by the assumed condition (i) or (ii) of this corollary. The stated result then follows from Theorem \ref{prop:quadraticgrow}.
\endproof

\subsection{The asymptotically superlinear convergence rate of the ALM}

Based on a recent paper \cite{CST2017}, the derived quadratic growth condition of \eqref{dual} guarantees the asymptotically superlinear convergence rate of the  KKT residual of the iterative sequence generated by the ALM in \eqref{ALM:iter} for solving the DNN projection problem under easy-to-implement stopping criteria.

In the seminal paper of Rockafellar \cite{rockafellar1976augmented}, he suggested the following stopping criteria for the inexact computation of the augmented Lagrangian subproblems:
\begin{eqnarray*}
(A)\quad  &&f_k(X^{\,k+1}) - \inf f_k \,\leq\, \varepsilon_k^2\,/\,2\sigma_k, \\[5pt]
(B) \quad && f_k(X^{\,k+1}) - \inf f_k \,\leq\, \left(\,\eta_k^2\,/\,2\sigma_k\,\right)\left\|\,\left(S^{\,k+1}-S^{\,k}, Z^{\,k+1}-Z^{\,k}\right)\,\right\|^2,
\end{eqnarray*}
where  $\{\varepsilon_k\}$ and $\{\eta_k\}$ are two positive summable sequences. In particular, the criterion $(A)$ is sufficient to ensure the global convergence of the dual variable sequence $\{(S^{\,k}, Z^{\,k})\}$ to a multiplier of \eqref{eq:DNN}, while the criterion $(B)$, together with the dual quadratic growth condition, ensure its asymptotic superlinear convergence rate. It may be difficult to execute $(A)$ and $(B)$ for general convex problems since the value $\inf\, f_k$ is generally unknown. One nice feature of the augmented Lagrangian subproblem in \eqref{ALM:iter} is that the function $f_k(\cdot)$ is continuously differentiable and strongly convex with modulus $1$ for any $k\geq 0$. Therefore, it holds that
$$
f_k(X) - \inf f_k \,\leq\, \frac{1}{2}\,\left\|\,\nabla f_k(X)\,\right\|^2,\quad \forall \; X\in \mathbb{S}^n.
$$
The above inequality is adopted from~\cite[(4.5)]{rockafellar1976augmented}, which has its source from the proof of \cite[Proposition 2]{Kort1976}. As a consequence, the criteria $(A)$ and $(B)$ can be executed by
\begin{eqnarray*}
(A')\quad  &&\|\nabla f_k({X^{k+1}})\| \,\leq\, \varepsilon_k\,/\sqrt{\sigma_k}, \\[5pt]
(B') \quad && \|\nabla f_k({X^{k+1}})\| \,\leq\, \left(\,\eta_k\,/\sqrt{\sigma_k}\,\right)\left\|\,\left(S^{\,k+1}-S^{\,k}, Z^{\,k+1}-Z^{\,k}\right)\,\right\|.
\end{eqnarray*}

The following theorem states the global convergence and the asymptotically superlinear convergence rate of the ALM for solving \eqref{eq:DNN} under criteria $(A')$ and $(B')$.

\begin{theorem}\label{thm:alm}
	Let $\left\{\left(\,X^{\,k}, S^{\,k},Z^{\,k}\,\right)\right\}$ be an infinite  sequence generated by the ALM   with stopping criterion $(A')$. Then the whole sequence $\left\{\left(\,X^{\,k}, S^{\,k},Z^{\,k}\,\right)\right\}$ is bounded with $\{X^{\,k}\}$ converging to the unique primal optimal solution $X^\infty$ and $\left\{\left(\,S^{\,k},Z^{\,k}\,\right)\right\}$ converging to some point $(S^\infty, Z^\infty)\in {\rm SOL}_{\bf D}$.
	
	If the criterion $(B')$ is also executed and the dual quadratic growth condition	holds at $(S^\infty, Z^\infty)$  with modulus $\kappa$, then there exists $k_0 \geq 0$ such that  for all $k\geq k_0$, 
	\begin{subequations}
		\begin{align}
		{\rm dist}\left(\left(S^{\,k+1},Z^{\,k+1}\right),  {\rm SOL}_{\rm D}\right)\,\leq\, \mu_k\,{\rm dist}\left(\left(S^{\,k},Z^{\,k}\right), {\rm SOL}_{\rm D}\right),\label{rate:dual}\\[3pt]
		\left\|\,\Res\left(X^{\,k+1}, S^{\,k+1}, Z^{\,k+1}\right)\,\right\|\,\leq \, \mu_k'\,{\rm dist}\left(\left(S^{\,k},Z^{\,k}\right), {\rm SOL}_{\rm D}\right),\label{rate:kkt}
		\end{align}
	\end{subequations}
	where the function $\Res(\,\cdot\,)$ is defined in \eqref{defn:kkt} and the constants $\mu_k, \mu_k'$ are given by
	$$
	\left\{\begin{array}{ll}
	\mu_k\,\triangleq \, \left[\,{\eta}_k+({\eta}_k+1)/\sqrt{1+ \sigma_k^2\,\kappa^2}\,\right]/(1-{\eta}_k)\to\mu_\infty \,\triangleq \, 1/\sqrt{1+ \sigma_\infty^2\,\kappa^2}\,,\\[0.1in]
	\mu_k' \,\triangleq \,  [\,\eta_k/\sqrt{\sigma_k}+2/\sigma_k\,]/(1-{\eta}_k)\to \mu'_\infty \,\triangleq \, 2/\sigma_\infty\,.
	\end{array}\right.
	$$
	Moreover, $\mu_\infty = \mu_\infty' = 0$ if $\sigma_{\infty} = +\infty$.
\end{theorem}
\proof
Since the Slater condition of problem \eqref{eq:DNN} trivially holds, the solution set of the dual problem is nonempty. Then the global convergence of $\{X^{\,k}\}$ and $\{(S^{\,k},Z^{\,k})\}$ follows from \cite[Theorem 4]{rockafellar1976augmented}. The inequality \eqref{rate:dual} under criterion $(B')$ is due to \cite[Theorem 4]{rockafellar1976augmented}. The inequality \eqref{rate:kkt} can be obtained from \cite[Theorem 2]{CST2017}.
\endproof

The above theorem shows that under the dual quadratic growth condition, the dual sequence generated by the ALM converges Q-linearly and the KKT residual of the primal-dual sequence converges R-linearly if $\displaystyle\lim_{k\to\infty} \sigma_k < +\infty$. The linear convergence rates $\mu_k$ and $\mu_k'$ can be  arbitrarily small with a sufficiently large value of $\sigma_k$. This type of convergence property is called ``arbitrarily fast linear convergence'' by Powell in \cite{Powell1972} when he studied the ALM for solving equality constrained nonlinear programming. The convergence rate of the dual sequence becomes asymptotically superlinear when $\sigma_k\to +\infty$. It is this property that distinguishes the ALM from various first order methods such as the  alternating direction method of multipliers (ADMM), where the latters' linear convergence rate (established under primal-dual type error bound conditions) is always close to $1$ for ill-conditioned problems; see, e.g., \cite[Theorem 2]{han2015lineara} for the convergence rate of the ADMM.

\section{A Semismooth Newton-CG Based Augmented Lagrangian Method for ({\bf P})}
\label{sec:SSN}

In this section, we discuss the semismooth Newton-CG method for solving the augmented Lagrangian subproblems in \eqref{ALM:iter}.

Recall that a  locally Lipschitz continuous function $F: \mathcal{O} \subseteq X \to Y$ defined on an open set $\mathcal{O}$ is said to be semismooth at $x \in  \mathcal{O}$ if $F$ is directionally differentiable at $x$ and for any $V \in \partial F(x+\Delta x)$ with $\Delta x\to 0$,
\[
F (x + \Delta x) - F (x) - V \Delta x = o(\|\Delta x\|),
\]
and $F$ is said to be strongly semismooth at $x$ if $F$ is semismooth at $x$ and 
\[
F(x + \Delta x) - F(x) - V \Delta x = O(\|\Delta x\|^2).
\]
$F$ is said to be a semismooth (respectively, strongly semismooth) function on $\mathcal{O}$ if it is semismooth (respectively, strongly semismooth) everywhere in $\mathcal{O}$.

Given a positive penalty parameter $\sigma$ and the dual variables $(\widetilde{S}, \widetilde{Z})\in \mathbb{S}^n\times \mathbb{S}^n$, the augmented Lagrangian subproblem is given by:
\begin{equation}\label{subALM}
\operatornamewithlimits{minimize}_{x\in \mathbb{S}^n}\, \tilde{f}(X)\,\triangleq\, L_{\sigma}(\,X;\, \widetilde{S},\widetilde{Z}\,).
\end{equation}
The function $\tilde{f}$ is continuously differentiable with the gradient given by
\begin{eqnarray}\label{eq:gradient}
\grad \tilde{f}(X) = X-G -  \Pi_{\mathbb{S}_+^n}(\widetilde{S} - \sigma X) -  \Pi_{\mathbb{N}^n}(\widetilde{Z} - \sigma X), \quad X\in \mathbb{S}^n.
\end{eqnarray}
Thus,  the optimal solution to the subproblem \eqref{subALM} can be obtained via the solution of the nonlinear equation
$$
\nabla \tilde{f}(X) = 0, \quad X\in \mathbb{S}^n.
$$
It is known from \cite{Sun2002semismooth} and \cite[Proposition 7.4.4 \& 7.4.7]{facchinei2007finite} that $\nabla \tilde{f}:\mathbb{S}^n\to\mathbb{S}^n$ is strongly semismooth so that the semismooth Newton-CG method is applicable to solve the above equation. The generalized Jacobian of $\nabla \tilde{f}$ at $X  \in \mathbb{S}^n$ is given by
$$
\partial (\nabla \tilde{f}) (X)  = \Big\{
I +\sigma ( V_1 + V_2) \mid V_1 \in \partial  \Pi_{\mathbb{S}_+^n}(\widetilde{S} - \sigma X),\; V_2 \in \partial  \Pi_{\mathbb{N}^n}(\widetilde{Z} - \sigma X)
\Big\}.
$$
A globally convergence semismooth Newton-CG method with line search is described as follows. 
\pagebreak

\noindent\makebox[\linewidth]{\rule{\textwidth}{1pt}}

\noindent{\bf Algorithm  SNCG:} a {\bf S}emismooth {\bf N}ewton-{\bf CG} method for solving the subproblem of the ALM

\noindent\makebox[\linewidth]{\rule{\textwidth}{1pt}}

\noindent
{\bf Initialization.} Given $\mu \in (0,1/2)$, $\eta\in (0,1)$, $\tau \in (0,1]$ and $\delta  \in (0,1)$. Iterate the
following steps for $j\geq 0$. \\[5pt]
{\bf Step 1.} Choose $V_1^j \in \partial  \Pi_{\mathbb{S}_+^n}(\widetilde{S} - \sigma X)$ and $V_2^j \in \partial  \Pi_{\mathbb{N}^n}(\widetilde{Z} - \sigma X)$. Solve the following linear system to find $\Delta X_j$ by the conjugate gradient method:
\[
(I + \sigma V_1^j + \sigma V_2^j)\Delta X+ \nabla \tilde{f} (X^j) = 0,
\]
until
$$
\|(I + \sigma V_1^j + \sigma V_2^j)\Delta X_j +\nabla \tilde{f}(X^j)\|\leq \min\,\left(\,\eta, \, \|\nabla \tilde{f}(X^j)\|^{1+\tau}\,\right).
$$
{\bf Step 2.} (Line search)
Set $\alpha_j = \delta^{m_j}$, where $m_j$ is the first nonnegative integer $m$ for which
$$
\tilde{f}(X^j +\delta^m \Delta X_j) \leq \tilde{f}(X^j) + \mu \delta^m\langle \nabla \tilde{f}(X^j), \Delta X_j\rangle.
$$
{\bf Step 3.}
Set $X^{j+1} = X^j + \alpha_j \Delta X_j$.

\noindent\makebox[\linewidth]{\rule{\textwidth}{1pt}}

The global convergence and the superlinear convergence rate are stated in the following proposition, whose proof can be established similarly as in \cite[Theorems 3.4 and 3.5]{zhao2010newton}.
\begin{proposition}
	Let the sequence $\{X^j\}$ be generated by Algorithm SNCG. Then $\{X^j\}$ converges to the unique optimal solution $\overline{X}$ of the problem in \eqref{subALM} and
	$$
	\| X^{j+1} - \overline{X}\| = O(\|X^j -\overline{X}\|^{ 1+\tau}).
	$$
\end{proposition}

The major computational cost of the SNCG method is to solve the following linear system
\begin{equation}\label{eq:newton}
(I + \sig V_1 + \sig V_2) (\Delta X) = R
\end{equation}
by the conjugate gradient method,
where $R\in\S^n$ is a given right-hand-side. In the following, we provide  particular choices of $V_1 \in \partial  \Pi_{\mathbb{S}_+^n}(\widetilde{S} - \sigma X)$ and  $V_2 \in \partial  \Pi_{\mathbb{N}^n}(\widetilde{Z} - \sigma X)$ with explicitly expressions of the products $V_1(\Delta X)$ and $V_2(\Delta X)$ with any $\Delta X\in \mathbb{S}^n$. Let $\lambda_1\geq \lambda_2\geq \ldots\geq \lambda_n$ be the eigenvalues of $\widetilde{S} - \sigma X$ and
$P$ be a corresponding orthogonal matrix of eigenvectors, i.e.,
$$
\widetilde{S} - \sigma X = P{\rm diag}(\lambda_1, \ldots, \lambda_n)P^T.
$$
We also denote the following index sets
$$
\alpha\,\triangleq \,\{i\mid \lambda_i>0, \; 1\leq i\leq n\}\quad \mbox{and}\quad \bar{\alpha} = \{i\mid \lambda_i \leq 0, \; 1\leq i\leq n\}.
$$
Denote the matrix $\Omega\in \S^n$ as
$$
\Omega \,\triangleq \, \left(\begin{array}{cc}
E_{\alpha\alpha}& \nu_{\alpha\bar{\alpha}} \\[2pt]
\nu^{\,T}_{\alpha\bar{\alpha}} & 0
\end{array}
\right)\quad {\rm with}\quad (E_{\alpha\alpha})_{ij} \,\triangleq\, 1,\; i, j\in \alpha\;\,\mbox{and}\;\, (\nu_{\alpha\bar\alpha})_{ij}\,\triangleq\, \frac{\lambda_i}{\lambda_i - \lambda_j},\; i\in \alpha, \; j \in\bar{\alpha}.
$$
In addition, we write the matrix $M\in \S^n$ as
$$
M_{ij} = \left\{\begin{array}{ll}
1 &  \quad {\rm if} \;\, (\widetilde{Z} - \sigma X)_{ij}\geq 0, \\[3pt]
0 & \quad  {\rm otherwise},
\end{array}
\right. i,j=1,\ldots,n.
$$
Based on the above preparations,  the linear operators $V_1$ and $V_2$ are chosen such that
\begin{eqnarray*}
	V_1(\Delta X) = P \left[\Omega\circ(\,P^{\,T} \Delta X \,  P\,)\right]P^{\,T}, \quad
	V_2(\Delta X) = M\circ \Delta X, \quad  \Delta X\in\S^n,
\end{eqnarray*}
where ``$\circ$" denotes the Hadamard product between two matrices. Moreover, if we partition $P$ corresponding to $\alpha$ and $\bar{\alpha}$, namely $ P = [P_1,P_2] $, then by making use of the special structure of $\Omega$, we have that
\begin{align*}
V_1(\Delta X) &= \begin{bmatrix}
P_1 &P_2
\end{bmatrix}\left( \begin{bmatrix}
E_{\alpha\alpha} & \nu_{\alpha\bar{\alpha}} \\ \nu_{\alpha\bar{\alpha}}^T & 0
\end{bmatrix}\circ \begin{bmatrix}
P_1^T\Delta XP_1 & P_1^T\Delta X P_2 \\ P_2^T \Delta XP_1 & P_2^T\Delta XP_2^T
\end{bmatrix} \right)\begin{bmatrix}
P_1^T \\ P_2^T
\end{bmatrix}
\\[5pt]
& = P_1P_1^T\Delta XP_1P_1^T + P_1[\nu_{\alpha\bar{\alpha}}\circ (P_1^T\Delta X P_2)]P_2^T + P_2[\nu_{\alpha\bar{\alpha}}^T\circ( P_2^T \Delta XP_1)]P_1^T.
\end{align*}

In order to reduce the iteration number of the conjugate gradient method for solving \eqref{eq:newton}, one may consider the following preconditioned system:
\[
[\,I + \sig(I+\sig V_1)^{-1} V_2\,] (\Delta X) = - (I+\sig V_1)^{-1} R.
\]
It can be shown that
\begin{eqnarray*}
	(I+\sig V_1)^{-1} (\Delta X) = P\, [\,\Xi\circ(P^{\,T} \Delta XP)\,]\,P^{\,T},
	\quad \Xi_{ij} = \frac{1}{1+\sigma \Omega_{ij}}, \; i,j=1,\ldots,n.
\end{eqnarray*}
Moreover, to make use of the $(2,2)$ block of zeros in $\Omega$ to reduce the computational cost, we may rewrite the above computation as
\begin{eqnarray*}
	(I+\sig V_1)^{-1} (\Delta X) = \Delta X -P[\,\Sigma \circ(P^{\,T} \Delta X P)\,]\, P^{\,T},
	\quad \Sigma_{ij} = \frac{\sig\Omega_{ij}}{1+\sigma \Omega_{ij}}, \; i,j=1,\ldots,n.
\end{eqnarray*}
On the other hand, if one wants to make use of the $(1,1)$ block of ones in $\Omega$ to reduce the computational cost, we may use the following computation:
\begin{eqnarray*}
	(I+\sig V_1)^{-1} (\Delta X) = \frac{1}{1+\sig}\Big(\Delta X + P [\, \Theta \circ(P^{\,T} \Delta X P)\,]\, P^{\,T}\Big),
	\quad \Theta_{ij} = \frac{\sig(1-\Omega_{ij})}{1+\sigma \Omega_{ij}}, \; i,j=1,\ldots,n.
\end{eqnarray*}

\section{Numerical Experiments}
\label{sec:numerical}
In this section, we conduct extensive numerical experiments to compare the performance of several methods on different data sets.

In our numerical experiments, we adopt an accelerated proximal gradient method (APG) of Nesterov \cite{Nesterov1983} to warm start the ALM. Let $(S^*, Z^*)\in\mathbb{S}^n\times \mathbb{S}^n$ be an optimal solution of problem (\ref{dual}). It is easy to derive that $(S^*, Z^*)$ always satisfies  $Z^* = \Pi_{\mathbb{N}^n}( - G - S^*)$ and
\begin{equation}\label{apg:opt}
S^* \in  \displaystyle\operatornamewithlimits{argmin}_{S\in \mathbb{S}^n}\, \left\{\phi(S)\,\triangleq \,\frac{1}{2}\|\Pi_{\mathbb{N}^n}(S + G)\|^2\mid S\in\mathbb{S}_+^n\right\}.
\end{equation}
Thus, we can eliminate the variable $Z$ in \eqref{dual} and solve the single-variable problem \eqref{apg:opt} in terms of $S$. Observe that the function $\phi$ is continuously differentiable with the gradient given by
$$
\nabla \phi(S) = \Pi_{\mathbb{N}^n}(S + G), \quad S\in\mathbb{S}^n.
$$
Moreover, the following inequality holds due to the global Lipschitz continuity (with modulus $1$) of the projection operator $\Pi_{\mathbb{N}^n}(\cdot)$:
$$
\phi(S)\, \leq \, \widehat{\phi}(S;\widehat{S})\,\triangleq \,\phi(\widehat{S}) + \left\langle \nabla \phi(\widehat{S}), S - \widehat{S}\right\rangle + \frac{1}{2}\,\|\,S - \widehat{S}\,\|^2, \quad\forall\; S, \widehat{S}\in\mathbb{S}^n.
$$
Given an initial point $\widetilde{S}^0 = S^0\in\mathbb{S}^n$ and parameter $t_0 = 1$, we use a variant of the accelerated proximal gradient method given in \cite{Beck2009} that executes the following iterative steps:
$$
\left\{\begin{array}{ll}
S^{\,k+1} =  \displaystyle\operatornamewithlimits{argmin}_{S\in \mathbb{S}^n}\, \left\{\widehat{\phi}(S; \widetilde{S}^{\,k})\mid S\in\mathbb{S}_+^n\right\} = \Pi_{\mathbb{S}_+^n}\left[\,\widetilde{S}^{\,k} - \Pi_{\mathbb{N}^n}(\widetilde{S}^{\,k}+G)\,\right],\\[0.15in]
t_{k+1} = \displaystyle\frac{1}{2}\left(\,1+\sqrt{1+4t_k^2}\,\right),\\[0.15in]
\widetilde{S}^{k+1} = S^{\,k+1}  + \displaystyle\frac{t_k-1}{t_{k+1}}\,(S^{\,k+1} - S^{\,k}).
\end{array}\right.
$$

For comparison purposes, we also test Dykstra's algorithm \cite{Dykstra1983} and the alternating direction method of multiplier (ADMM) to solve \eqref{eq:DNN}. Dykstra's algorithm \cite{Dykstra1983} is a  variant of the alternating projection method for computing the projection onto the intersection of a finite  number of closed convex sets. It is well known that Dykstra's algorithm is a particular block coordinate descent method applied to the dual problem \eqref{dual} \cite{Gaffke1989acyclic}. In order to apply the ADMM to \eqref{eq:DNN}, we first reformulated the problem as
$$\begin{array}{ll}
\min  & \displaystyle\frac{\alpha}{2}\,\|\,X_1 - G\,\|^2 + \frac{1-\alpha}{2}\,\|\,X_2-G\,\|^2\\[0.2in]
{\rm s.t.} & X_1 - X_2 = 0, \; X_1\in\mathbb{S}_+^n, \;X_2\in \mathbb{N}^n,
\end{array}
$$
where $\alpha\in (0,1)$ is a given parameter. Let $\sigma$ be a positive penalty parameter. The corresponding augmented Lagrangian function of the above problem is given by
$$
L_\sigma(X_1,X_2,W) = \displaystyle\frac{\alpha}{2}\,\|\,X_1 - G\,\|^2 + \frac{1-\alpha}{2}\,\|\,X_2-G\,\|^2 + \langle X_1 - X_2, W\rangle  + \frac{\sigma}{2}\|X_1 - X_2\|^2.
$$
Given initial points $X_2^0$ and  $W^0$ in $\mathbb{S}^n$ and a positive penalty parameter $\sigma$, the $(k+1)$-th iteration of the ADMM is given by
\[
\left\{\begin{array}{ll}
X_1^{k+1} = \displaystyle\operatornamewithlimits{argmin}_{X_1\in \mathbb{S}^n} L_\sigma(X_1, X_2^k,W^k) = \Pi_{\mathbb{S}_+^n}\left[(\alpha + \sigma)^{-1}\left(\alpha G + \sigma X_2^k - W^k\right)\right],\\[0.1in]
X_2^{k+1}  = \displaystyle\operatornamewithlimits{argmin}_{X_2\in \mathbb{S}^n} L_\sigma(X_1^{k+1}, X_2,W^{k}) = \Pi_{\mathbb{N}^n}\left[(1-\alpha + \sigma)^{-1}\left((1-\alpha) G + \sigma X_1^{k+1} + W^k\right)\right],\\[0.1in]
W^{k+1} = W^k + \tau\sigma\left(X_1^{k+1} - X_2^{k+1}\right),
\end{array}\right.
\]
where $\tau\in (0,\frac{\sqrt{5}+1}{2})$ is the step-length. We take $\tau = 1.618$ in our numerical experiments.

We terminate all the algorithms if the relative KKT residual
$$
\begin{array}{ll}
\eta \,\triangleq \, \displaystyle\frac{1}{\max\{1,\|G\|\}}\max\left\{\begin{array}{cc}
\|X^{\,k} - G - S^{\,k} - Z^{\,k}\|\\[3pt]
\|X^{\,k} - \Pi_{\mathbb{S}_+^n}(X^{\,k})\|,\;\|S^{\,k} - \Pi_{\mathbb{S}_+^n}(S^{\,k})\|,\;
|\langle X^{\,k},S^{\,k}\rangle|/(1+\|S^{\,k}\|)\\[3pt]
\|X^{\,k} - \Pi_{\mathbb{N}^n}(X^{\,k})\|,\;\|Z^{\,k} - \Pi_{\mathbb{N}^n}(Z^{\,k})\|,\;\,
|\langle X^{\,k},Z^{\,k}\rangle|/(1+\|Z^{\,k}\|)
\end{array}\right\}
\leq {\rm tol},
\end{array}
$$
where the tolerance ``tol'' is set to be $10^{-12}$ in the experiments. The algorithms will also be stopped when they reach the maximum number of iterations ($200$ for the ALM, and $20,000$ for the APG, the ADMM and Dykstra's algorithm).

In the rest of this section, we conduct experiments with input matrix $G$ generated from synthetic and real data. The number of iterations, the final KKT residuals and the computational time for each method are reported. For the  ALM, we also report the total number of semismooth Newton iterations needed to solve the ALM subproblems and the  number of APG iterations taken for the purpose of warm-starting. For instance, the item $50(257, 1190)$ in the first row under the column `alm' in Table \ref{tab-hankel} means that the number of ALM iterations is 50 with a total of $257$ semismooth iterations and $1190$ APG iterations. The computational time is in the format of ``hours:minutes:seconds''. We also check the strict complementarity with respect to the positive semidefinite constraint at the approximate KKT solution $ (\overline{X},\overline{S},\overline{Z})$ given by the last iterate of the ALM algorithm that is defined by the quantity
\[
\text{sc}:=\frac{\lambda_{\min}(\overline{X}+\overline{S})}{\lambda_{\max}(\overline{X}+\overline{S})},
\]
where $\lambda_{\min}(\overline{X}+\overline{S})$ and $\lambda_{\max}(\overline{X}+\overline{S})$ denote the minimal and maximal eigenvalues of $ \overline{X}+\overline{S} $, respectively.  If the quantity ``$\text{sc}$'' is substantially larger than \text{tol}, then one can confidently conclude that $\mbox{rank}(\overline{X}) +\mbox{rank}(\overline{S}) = n$, which implies that the quadratic growth condition holds at $\overline{X}$ due to Corollary~\ref{coro:sdpTphi}. However, it is worth  mentioning that in order for the dual quadratic growth condition to hold at $\left(\overline{S}, \overline{Z}\right)$, we only need the existence of a dual solution pair $\left(\widehat{S}, \widehat{Z}\right)$ such that  $\mbox{rank}(\overline{X}) +\mbox{rank}(\widehat{S}) = n$. Unfortunately, the latter condition is difficult to verify numerically.

All experiments are run in {\sc Matlab} R2018b on a workstation with Intel Xeon processor E5-2680v3 @2.50GHz (12 cores and 24 threads) and 128GB of RAM, equipped with 64-bit Windows 10 OS.

\subsection{Experiments on synthetic data} \label{subsection-random}
We first conduct experiments on four classes of synthetic data: matrices whose projections are zeros, Hankel matrices, randomly generated noisy low rank sparse matrices and Toeplitz matrices.

\vskip 0.1in
\noindent
{\bf Example 1: matrices whose projections have zero solutions}. For given $ S\in\S^n_+ $ (with $ \mathrm{rank}(S) < n $) and  nonnegative  $ Z\in \S^n $  (with $ \mathrm{rank}(Z) < n $), let $ G = -(S+Z) $. Obviously, $ (X=0,S,Z) $ satisfies $ \mathcal{R}(X,S,Z) = 0 $. In our experiments for Table \ref{tab-Zero}, we generate both matrices $ S $ and $ Z $ randomly via the following {\sc Matlab} script:
\begin{verbatim}
    Stmp = randn(n, 2); S = Stmp*Stmp';
    Ztmp = rand(n, 2);  Z = Ztmp*Ztmp';
    G = -(S+Z);
    G = G/norm(G,'fro');
\end{verbatim}

\noindent
{\bf Example 2: Hankel matrices}.
A Hankel matrix is a square matrix in which each ascending skew-diagonal form left to right is constant. In Table \ref{tab-hankel}, we consider Hankel matrices with dimension $ n $ generated by the following  {\sc Matlab} commands:
\begin{verbatim}
    G = hankel(-(1:n)',(1:n)');
    G = G/norm(G,'fro');
\end{verbatim}

\noindent
{\bf Example 3: noisy low rank sparse matrices}.
In our numerical experiments in Table \ref{tab-sparse}, the input matrices $G$ are
noisy low rank sparse matrices  generated via the following {\sc Matlab} commands:
\begin{verbatim}
    V = sprand(n,10,0.5);
    G0 = -V*V';  E = randn(n); E = 0.5*(E+E');
    G = 0.85*G0+0.15*E;
    G = G/norm(G,'fro');
\end{verbatim}

\noindent
{\bf Example 4: Toeplitz matrices}.
A Toeplitz matrix  is a matrix in which each descending diagonal from left to right is constant.
In our numerical experiments in Table \ref{tab-toeplitz}, the input matrices $G$ are Toeplitz matrices generated as follows:
\begin{verbatim}
    c = -rand(n,1);  c(1:n/25) = ones(n/25,1);
    G = toeplitz(c);
    G = G/norm(G,'fro');
\end{verbatim}

\begin{footnotesize}
	\begin{table}[ht!]
		\caption{Numerical results on synthetic data whose projection is the zero matrix.}
		\label{tab-Zero}
		\resizebox{\textwidth}{!}{
			\begin{tabular}{|c|c|c|c|c|}\hline
				& Iteration & KKT residual & time & sc \\ \hline
				n & alm $|$ apg $|$ admm $|$ dykstra & alm $|$ apg $|$ admm $|$ dykstra & alm $|$ apg $|$ admm $|$ dykstra & alm \\ \hline
				400 &    5(  25,  550) $|$ 1430 $|$ 20000 $|$ 20000 & 2.8e-14 $|$ 9.6e-13 $|$ 4.4e-10 $|$ 3.6e-07 &13$|$26$|$6:42$|$5:52 & 1.1e-15\\ \hline
				600 &   13(  70,  900) $|$ 3140 $|$ 20000 $|$ 20000 & 7.3e-14 $|$ 5.2e-13 $|$ 1.2e-08 $|$ 6.2e-07 &45$|$1:50$|$11:56$|$11:16 & 2.6e-15\\ \hline
				800 &   10(  50,  540) $|$ 2130 $|$ 20000 $|$ 20000 & 7.4e-13 $|$ 9.9e-13 $|$ 5.8e-09 $|$ 4.3e-07 &52$|$2:16$|$23:32$|$20:55 & 2.0e-14\\ \hline
				1000 &   11(  55,  760) $|$ 4940 $|$ 20000 $|$ 20000 & 9.2e-13 $|$ 9.4e-13 $|$ 9.9e-09 $|$ 5.0e-07 &1:43$|$8:02$|$35:43$|$30:53 & 1.6e-15\\ \hline
				1200 &   20( 100,  780) $|$ 7210 $|$ 20000 $|$ 20000 & 3.5e-14 $|$ 9.9e-13 $|$ 2.7e-08 $|$ 6.2e-07 &3:21$|$18:00$|$57:36$|$53:47 & 2.0e-15\\ \hline
				1400 &   34( 215, 1600) $|$ 9290 $|$ 20000 $|$ 20000 & 4.7e-13 $|$ 8.8e-13 $|$ 1.0e-08 $|$ 4.7e-07 &10:41$|$35:15$|$  1:24:11$|$  1:14:42 & 5.9e-16\\ \hline
		\end{tabular}}
	\end{table}
\end{footnotesize}

\begin{footnotesize}
	\begin{table}[H]
		\centering
		\caption{Numerical results on Hankel matrices.}
		\label{tab-hankel}
		\resizebox{\textwidth}{!}{
			\begin{tabular}{|c|c|c|c|c|}\hline
				& iteration & KKT residual & time & sc\\ \hline
				$n$ & alm $|$ apg $|$ admm $|$ dykstra & alm $|$ apg $|$ admm $|$ dykstra & alm $|$ apg $|$ admm $|$ dykstra &  alm \\ \hline
				400  &   50( 257, 1190) $|$ 20000 $|$ 20000 $|$ 20000 & 9.2e-13 $|$ 8.6e-12 $|$ 5.5e-11 $|$ 6.6e-08 &50$|$5:19$|$5:22$|$4:48                   & 5.0e-17 \\ \hline
				600  &   40( 202, 1050) $|$ 20000 $|$ 20000 $|$ 20000 & 9.4e-13 $|$ 7.5e-12 $|$ 6.5e-11 $|$ 6.7e-08 &1:36$|$11:00$|$11:21$|$10:32              & 2.5e-15 \\ \hline
				800  &   60( 302, 1010) $|$ 20000 $|$ 20000 $|$ 20000 & 5.7e-13 $|$ 6.7e-12 $|$ 1.2e-10 $|$ 6.7e-08 &3:53$|$19:11$|$20:14$|$18:14              & 2.8e-16 \\ \hline
				1000 &   45( 227, 1050) $|$ 20000 $|$ 20000 $|$ 20000 & 9.0e-13 $|$ 5.0e-12 $|$ 2.0e-10 $|$ 6.7e-08 &5:25$|$29:34$|$32:00$|$28:21              & 1.2e-16 \\ \hline
				1200 &   50( 252, 1090) $|$ 20000 $|$ 20000 $|$ 20000 & 9.9e-13 $|$ 7.5e-12 $|$ 2.4e-10 $|$ 6.8e-08 &8:58$|$45:27$|$48:07$|$43:57              & 1.4e-15 \\ \hline
				1400 &   65( 342, 1370) $|$ 20000 $|$ 20000 $|$ 20000 & 8.3e-13 $|$ 7.6e-12 $|$ 3.2e-10 $|$ 6.8e-08 &19:21$|$  1:04:34$|$  1:08:26$|$  1:02:47 & 4.0e-16 \\ \hline
		\end{tabular}}
	\end{table}
\end{footnotesize}

\begin{footnotesize}
	\begin{table}[H]	
		\centering
		\caption{Numerical results on noisy low rank sparse matrices.}
		\label{tab-sparse}
		\resizebox{\textwidth}{!}{
			\begin{tabular}{|c|c|c|c|c|}\hline
				& iteration & KKT residual & time & sc\\ \hline
				$n$ & alm $|$ apg $|$ admm $|$ dykstra & alm $|$ apg $|$ admm $|$ dykstra & alm $|$ apg $|$ admm $|$ dykstra & alm\\ \hline
				400  &  120( 603, 1320) $|$ 20000 $|$ 20000 $|$ 20000 & 5.8e-13 $|$ 3.3e-11 $|$ 2.9e-10 $|$ 2.2e-07 &1:50$|$6:22$|$6:45$|$5:51                 & 1.7e-09  \\ \hline
				600  &   45( 223,  980) $|$ 20000 $|$ 11661 $|$ 20000 & 9.4e-13 $|$ 1.7e-11 $|$ 8.7e-13 $|$ 7.4e-08 &1:56$|$11:53$|$7:22$|$10:55               & 2.1e-14  \\ \hline
				800  &   85( 428,  680) $|$ 20000 $|$ 18901 $|$ 20000 & 8.2e-13 $|$ 1.5e-11 $|$ 8.6e-13 $|$ 5.6e-08 &5:48$|$19:47$|$19:49$|$18:30              & 3.6e-14  \\ \hline
				1000 &   90( 453,  830) $|$ 20000 $|$ 20000 $|$ 20000 & 8.1e-13 $|$ 1.6e-11 $|$ 2.1e-11 $|$ 4.3e-08 &10:16$|$30:54$|$34:15$|$29:26             & 2.9e-16  \\ \hline
				1200 &   95( 478,  880) $|$ 20000 $|$ 20000 $|$ 20000 & 9.5e-13 $|$ 1.8e-11 $|$ 2.8e-11 $|$ 4.8e-08 &16:58$|$49:48$|$52:47$|$45:27             & 2.3e-15  \\ \hline
				1400 &   70( 353,  728) $|$ 20000 $|$ 18261 $|$ 20000 & 4.9e-13 $|$ 2.4e-11 $|$ 9.8e-13 $|$ 5.7e-08 &18:33$|$  1:11:20$|$  1:07:39$|$  1:05:26 & 8.8e-15  \\
				\hline
		\end{tabular}}
	\end{table}
\end{footnotesize}

\begin{footnotesize}
	\begin{table}[H]
		\centering
		\caption{Numerical results on Toeplitz matrices.}
		\label{tab-toeplitz}
		\resizebox{\textwidth}{!}{
			\begin{tabular}{|c|c|c|c|c|}\hline
				& iteration & KKT residual & time & sc \\ \hline
				$n$ & alm $|$ apg $|$ admm $|$ dykstra & alm $|$ apg $|$ admm $|$ dykstra & alm $|$ apg $|$ admm $|$ dykstra & alm \\ \hline
				400 &   21( 202,  740) $|$ 20000 $|$ 4501 $|$ 20000 & 7.8e-13 $|$ 2.2e-11 $|$ 9.5e-13 $|$ 5.4e-09 &33$|$5:47$|$1:25$|$5:39 & 5.9e-06\\ \hline
				600 &   14( 132,  630) $|$ 17670 $|$ 4121 $|$ 20000 & 9.4e-13 $|$ 9.3e-13 $|$ 7.5e-13 $|$ 2.5e-09 &48$|$9:29$|$2:28$|$10:28 & 5.9e-06\\ \hline
				800 &   19( 181,  580) $|$ 20000 $|$ 4321 $|$ 20000 & 7.9e-13 $|$ 8.9e-12 $|$ 7.1e-13 $|$ 2.6e-09 &1:44$|$19:10$|$4:31$|$18:42 & 9.2e-06\\ \hline
				1000 &   16( 154,  660) $|$ 15270 $|$ 3281 $|$ 20000 & 9.3e-13 $|$ 7.5e-13 $|$ 8.6e-13 $|$ 6.4e-11 &2:39$|$22:36$|$5:31$|$29:07 & 2.6e-06\\ \hline
				1200 &   14( 133,  650) $|$ 15300 $|$ 2921 $|$ 20000 & 9.4e-13 $|$ 9.2e-13 $|$ 6.9e-13 $|$ 1.1e-11 &3:39$|$36:22$|$7:44$|$46:40 & 3.6e-07\\ \hline
				1400 &   16( 153,  580) $|$ 17820 $|$ 2961 $|$ 20000 & 8.9e-13 $|$ 9.9e-13 $|$ 6.3e-13 $|$ 6.3e-12 &5:28$|$  1:00:39$|$10:45$|$  1:07:35 & 1.9e-05\\ \hline
		\end{tabular}}
	\end{table}
\end{footnotesize}

Table~\ref{tab-Zero} presents the computational results for matrices whose projections are zero matrices. One can observe that our ALM is obviously more efficient and robust than other methods in terms of the KKT residual and computational time. Moreover, the APG solves all the instances successfully while both the ADMM and Dykstra's algorithm cannot solve these problems to the desired accuracy. However, the ADMM is better than the Dykstra's algorithm in terms of the KKT residual but shares similar performance in terms of computational time.

From Table~\ref{tab-hankel}, one can observe that our ALM outperforms other methods in terms of number of iterations, KKT residuals and computational time. In fact, our ALM is always able to return a highly accurate solution with much shorter computational time. All the other algorithms cannot successfully solve the instances within $ 20000 $ iterations, with the APG performing slightly better than the ADMM, and much better than Dykstra's algorithm.

For the results of noisy low rank sparse matrices reported in Table~\ref{tab-sparse},  our ALM again outperforms all the other three methods. Both the APG and Dykstra's algorithm fail to reach the desired accuracy for all the tested instances within $20000$ iterations, while the ADMM can solve around half of the instances to the desired accuracy with much longer computational time compared to the ALM.

For the numerical results on Toeplitz matrices reported in  Table~\ref{tab-toeplitz}, one can find that the ALM still outperforms the other methods, but the ADMM also performs fairly well since it is only  about 2-3 times slower than the ALM in solving all the instances to the desired accuracy.

Notice that the term $\text{sc} \approx \text{tol} $ in Table~\ref{tab-Zero}, Table~\ref{tab-hankel} and Table~\ref{tab-sparse} but $\text{sc} \gg \text{tol}$ in Table~\ref{tab-toeplitz}. Even though we cannot conclude that the dual quadratic growth condition does not hold for the instances in the former three tables, we have indeed observed that the convergence rates in Table~\ref{tab-Zero} and Table~\ref{tab-hankel} are slower compared to those in Table~\ref{tab-toeplitz} in the sense that more ALM iterations are needed to solve the problems to the desired accuracy. Furthermore, from Table~\ref{tab-Zero}--Table~\ref{tab-sparse}, even though $ \text{sc} \approx \text{tol} $, the convergence rates in Table~\ref{tab-Zero} are faster compared to those in Table~\ref{tab-hankel} and Table~\ref{tab-sparse}. This again indicates that the failure of the strict complementarity condition at a particular point does not necessarily imply the failure of the quadratic growth condition. Therefore, as mentioned in the Introduction, one can see that the quadratic growth condition is quite mild.

\subsection{Experiments on DNN projection instances arising from solving the Lagrangian-DNN relaxations of quadratic optimization problems}

Problem \eqref{eq:DNN} arises naturally as a subroutine in the Lagrangian-DNN relaxation method for approximately solving a quadratic optimization problem (QOP) of the following form:
\begin{equation}\label{eq-QOP}
	\min_{\bm u}\;\;\left\{ \bm{u}^T\bm{Q}\bm{u}  + 2\bm{c}^T\bm{u} \;\middle|\; \begin{aligned}
	&\bm{u}\in \R_+^m,\;\bm{A}\bm{u}+\bm{b} = 0,\\
	&u_iu_j = 0 \;((i,j)\in \mathcal{E})
	\end{aligned} \right\},
\end{equation}
where $\bm{A}\in \R^{q\times m}$, $\bm{b}\in \R^q$, $ \bm{c}\in \R^m $, $ \bm{Q}\in \S^m $, and $ \mathcal{E}\subset \left\{ (i,j)\;:\;1\leq i<j\leq m \right\} $ are given data.
Let  $n = 1+m$, and
\[
\bm{Q}_0 := \begin{pmatrix}
0 & \bm{c}^T \\ \bm{c} & \bm{Q}
\end{pmatrix}\in \S^n,\; \bm{H}_0 := \begin{pmatrix}
1 & \bm{0}^T \\ \bm{0} & \bm{O}
\end{pmatrix}\in \S^n,\; \bm{Q}_{01} := \begin{pmatrix}
\bm{b}^T\bm{b} & \bm{b}^T\bm{A} \\ \bm{A}^T\bm{b} & \bm{A}^T\bm{A}
\end{pmatrix}\in \S^n,\;
\bm{Q}_{ij}  := \begin{pmatrix}
0 & \bm{0}^T \\ \bm{0} & \bm{C}_{ij}+\bm{C}_{ij}
\end{pmatrix}
\]
with $\bm{C}_{ij}$ being the $m\times m$ matrix whose $(i,j)$-th component is $1/2$ if $(i,j)\in \mathcal{E}$ and $ 0 $ otherwise.

It has been shown in \cite{Kim2016Lagrangian} that the nonconvex problem \eqref{eq-QOP} can be reformulated
as the following completely positive cone convex programming problem:
\[
\inf_{\bm X}\;\;\left\{\bm{Q}_0\bullet\bm{X}\;\middle|\; \bm{H}_0\bullet\bm{X} = 1,\; \bm{H}_1\bullet\bm{X} = 0,\; \bm{X}\in \mathbb{C}^{n,\ast}    \right\}
\]
where $\bm{H}_1:=\bm{Q}_{01}+\sum_{(i,j)\in \mathcal{E}}\bm{Q}_{ij}$ and
$X \bullet Y = {\rm tr}(XY)$ for any $X,Y\in\S^n$.
While the above problem is convex,  the conic constraint is unfortunately not computationally tractable.
As suggested in \cite{Kim2016Lagrangian}, one can approximately solve it via the following linearly constrained DNN relaxation problem
based on the fact that $\mathbb{C}^{n,\ast}  \subset \mathbb{D}^n$:
\begin{eqnarray}
\label{eq-QOP-DNN-P}
\inf_{\bm X}\;\;\left\{\bm{Q}_0\bullet\bm{X}\;\middle|\; \bm{H}_0\bullet\bm{X} = 1,\; \bm{H}_1\bullet\bm{X} = 0,\; \bm{X}\in \mathbb{D}^n    \right\}.
\end{eqnarray}
Its corresponding Lagrange dual problem  is given by
\begin{equation}\label{eq-QOP-DNN-D}
\sup_{y_0}\;\;\left\{y_0\;\middle|\; \bm{Z}+y_0\bm{H}_0+y_1\bm{H}_1 = \bm{Q}_0,\;\bm{Z}\in \mathbb{D}^{n,*},\; \bm{y} = (y_0,y_1)\in \R^2 \right\}.
\end{equation}
For the sake of computing a
lower bound of \eqref{eq-QOP} efficiently,
the authors in \cite{Kim2016Lagrangian} further considered the Lagrangian-DNN relaxations of \eqref{eq-QOP-DNN-P} and
its dual that are given by
\begin{align}
&\inf_{\bm X}\;\;\left\{ \bm{Q}_0\bullet\bm{X}+\lambda \bm{H}_1\bullet\bm{X}\;\middle|\; \bm{H}_0\bullet\bm{X} = 1,\;\bm{X}\in \mathbb{D}^n \right\} \label{eq-LagDNN-P}
\\
&\sup_{y_0}\;\;\left\{ y_0\;\middle|\; \bm{Q}_0+\lambda \bm{H}_1-y_0\bm{H}_0\in \mathbb{D}^{n,*}\right\}\label{eq-LagDNN-D},
\end{align}
where $\lambda > 0$ is a given Lagrangian parameter that should be chosen large enough to obtain a high quality relaxation. When  the bisection method is applied to solve
the problem \eqref{eq-LagDNN-D} for a given (large) $\lambda$,  the key step in each bisection iteration is to compute the following DNN projection
for any given $y \in \mathbb{R}$:
\[
\Pi_{\mathbb{D}^{n,*}}(\bm{G}_\lambda (y))= \bm{G}_\lambda (y) + \Pi_{\mathbb{D}^n}(-\bm{G}_\lambda (y))\quad \text{with}\;\; \bm{G}_\lambda (y):= \bm{Q}_0+\lambda \bm{H}_1-y\bm{H}_0.
\]
Therefore, an efficient solver for computing the projection onto the DNN cone is critical for solving the Lagrangian-DNN relaxation problem \eqref{eq-LagDNN-D}.

We conduct numerical experiments on DNN projection instances arising from the Lagrangian-DNN relaxation method for solving quadratic optimization problems \eqref{eq-QOP} associated with binary integer quadratic problems (BIQ) and  quadratic assignment problems (QAP). We set $ \lambda = 10^6\times\frac{\norm{\bm{Q}_0}}{\max(1,\norm{\bm{H}_1})} $ for all the experiments. The parameter $y$ is  chosen from the interval $[y^*\times \frac{1}{1000},y^*\times1000] $, where $y^*$ is the optimal solution for the dual conic relaxation problem \eqref{eq-QOP-DNN-D} that is known from the literature (see e.g., \cite{Kim2016Lagrangian}). Given $y$ and $\lambda$, we compute the matrix $\bm{G}_\lambda(y)$ and take its normalization (by the Frobenius norm) as the input matrix $G$. The test instances for the BIQ  and QAP problems are downloaded from BIQMAC library~(available at \url{http://www.biqmac.uni-klu.ac.at/biqmaclib.html}.) and QAPLIB~(available at \url{http://www.seas.upenn.edu/qaplib}.), respectively.

Tables~\ref{tab-bqp}--\ref{tab-tai} present the numerical results for all the four algorithms. It can be seen that the ALM is about two times faster than the APG in terms of the computational time when both of them reach the required accuracy level.  Compared with the APG, the ADMM solves half of the BIQ instances with much longer computational time and all the QAP instances with roughly the same efficiency. Dykstra's algorithm, on the other hand, cannot solve a large proportion of the  problems to the desired accuracy within $20000$ iterations. One also observes that for most of the QAP instances, it holds that $ \text{sc} \approx \sqrt{\text{tol}} $. This fact indicates that the strict complementarity condition is likely to hold for these problems and hence resulting in fast convergence rates as shown in the tables. This also shows that a strict complementarity solution is likely to exist in some real-world problems.

\begin{footnotesize}
	\begin{table}[H]
		\caption{Numerical results on bqp-data for BIQ problems. In this table, $y = -100$.}
		\label{tab-bqp}
		\resizebox{\textwidth}{!}{
			\begin{tabular}{|c|c|c|c|c|c|}\hline
				\multicolumn{2}{|c|}{}& iteration & KKT residual & time & sc\\ \hline
				problem & $n$ & alm $|$ apg $|$ admm $|$ dykstra & alm $|$ apg $|$ admm $|$ dykstra & alm $|$ apg $|$ admm $|$ dykstra & alm\\ \hline
				
				bqp250-2 & 	501 &   26( 127,  760) $|$ 2510 $|$ 12501 $|$ 20000 & 2.8e-13 $|$ 9.8e-13 $|$ 5.8e-13 $|$ 2.6e-06 &36$|$1:08$|$6:04$|$8:32       &6.7e-11   \\ \hline
				bqp250-4 & 	501 &   16(  79,  760) $|$ 2490 $|$ 12701 $|$ 20000 & 5.2e-13 $|$ 9.9e-13 $|$ 9.4e-13 $|$ 2.6e-06 &30$|$1:07$|$6:14$|$8:39       &5.4e-12   \\ \hline  		
				bqp250-6 & 	501 &   19(  92,  660) $|$ 2479 $|$ 12561 $|$ 20000 & 9.4e-13 $|$ 9.9e-13 $|$ 8.4e-13 $|$ 2.6e-06 &29$|$1:08$|$6:07$|$8:34       &1.1e-10   \\ \hline  		
				bqp250-8 & 	501 &   20(  97,  750) $|$ 2500 $|$ 12801 $|$ 20000 & 1.0e-12 $|$ 1.0e-12 $|$ 9.2e-13 $|$ 2.6e-06 &32$|$1:09$|$6:14$|$8:33       &1.0e-11   \\ \hline  	
				bqp250-10 & 	501 &   21( 103,  760) $|$ 2470 $|$ 12561 $|$ 20000 & 3.6e-13 $|$ 1.0e-12 $|$ 6.7e-13 $|$ 2.6e-06 &34$|$1:06$|$6:04$|$8:33      &4.4e-11\\ \hline  		
				bqp500-2 & 	1001 &   21( 101,  860) $|$ 2780 $|$ 20000 $|$ 20000 & 2.2e-13 $|$ 1.0e-12 $|$ 9.6e-09 $|$ 4.8e-06 &2:38$|$5:01$|$39:27$|$36:39  &1.5e-11   \\ \hline 	
				bqp500-4 & 	1001 &   16(  77,  900) $|$ 2659 $|$ 20000 $|$ 20000 & 9.5e-13 $|$ 9.8e-13 $|$ 9.6e-09 $|$ 4.8e-06 &2:28$|$4:49$|$40:02$|$36:35  &1.4e-11   \\ \hline  	
				bqp500-6 & 	1001 &   16(  78,  940) $|$ 2700 $|$ 20000 $|$ 20000 & 7.2e-13 $|$ 1.0e-12 $|$ 9.6e-09 $|$ 4.8e-06 &2:35$|$4:54$|$39:59$|$36:22  &6.1e-17   \\ \hline  	
				bqp500-8 & 	1001 &   16(  77,  940) $|$ 2680 $|$ 20000 $|$ 20000 & 3.9e-13 $|$ 1.0e-12 $|$ 9.6e-09 $|$ 4.8e-06 &2:36$|$4:52$|$39:57$|$36:41  &1.7e-11   \\ \hline		
				bqp500-10 & 	1001 &   15(  72,  910) $|$ 2710 $|$ 20000 $|$ 20000 & 8.8e-13 $|$ 9.9e-13 $|$ 9.6e-09 $|$ 4.8e-06 &2:32$|$5:00$|$39:56$|$36:56 &1.3e-11\\ \hline
		\end{tabular}}
	\end{table}
\end{footnotesize}

\begin{footnotesize}
	\begin{table}[H]
		\caption{Numerical results on bur-data for QAP problems. In this table, $ y = 10^4 $.}
		\label{tab-bur}
		\resizebox{\textwidth}{!}{
			\begin{tabular}{|c|c|c|c|c|c|}\hline
				\multicolumn{2}{|c|}{}	& iteration & KKT residual & time & sc \\ \hline
				problem & $n$ & alm $|$ apg $|$ admm $|$ dykstra & alm $|$ apg $|$ admm $|$ dykstra & alm $|$ apg $|$ admm $|$ dykstra & alm\\ \hline
				bur26a &	677 &   20(  38,  348) $|$ 1317 $|$ 1761 $|$ 20000 & 9.4e-13 $|$ 1.0e-12 $|$ 8.2e-13 $|$ 3.7e-12 &24$|$1:06$|$1:15$|$19:53  & 1.7e-06 \\ \hline
				bur26b &	677 &   19(  36,  310) $|$ 1198 $|$ 1821 $|$ 20000 & 8.6e-13 $|$ 1.0e-12 $|$ 6.8e-13 $|$ 5.0e-12 &22$|$58$|$1:16$|$19:28    & 2.5e-06 \\ \hline
				bur26e &	677 &   20(  39,  348) $|$ 1345 $|$ 1821 $|$ 20000 & 9.5e-13 $|$ 1.0e-12 $|$ 7.1e-13 $|$ 4.4e-12 &25$|$1:08$|$1:18$|$20:18  & 1.8e-06 \\ \hline
				bur26f &	677 &   19(  36,  310) $|$ 1210 $|$ 1761 $|$ 20000 & 8.4e-13 $|$ 1.0e-12 $|$ 8.9e-13 $|$ 5.8e-12 &23$|$1:02$|$1:15$|$19:43  & 2.6e-06 \\ \hline
				bur26g &	677 &   20(  40,  390) $|$ 1404 $|$ 1861 $|$ 20000 & 9.8e-13 $|$ 1.0e-12 $|$ 7.0e-13 $|$ 2.0e-12 &25$|$1:10$|$1:19$|$19:48  & 8.3e-07 \\ \hline
				bur26h &	677 &   21(  50,  449) $|$ 1385 $|$ 1821 $|$ 20000 & 9.5e-13 $|$ 1.0e-12 $|$ 5.8e-13 $|$ 2.8e-12 &30$|$1:09$|$1:19$|$19:40  & 1.2e-06 \\ \hline
		\end{tabular}}
	\end{table}
\end{footnotesize}

\begin{footnotesize}
	\begin{table}[H]
		\caption{Numerical results on chr-data for QAP problems. In this table, $ y = 10^5 $.}
		\label{tab-chr}
		\resizebox{\textwidth}{!}{
			\begin{tabular}{|c|c|c|c|c|c|}\hline
				\multicolumn{2}{|c|}{}	& iteration & KKT residual & time & sc\\ \hline
				problem &	$n$ & alm $|$ apg $|$ admm $|$ dykstra & alm $|$ apg $|$ admm $|$ dykstra & alm $|$ apg $|$ admm $|$ dykstra & alm\\ \hline
				chr20a &	401 &   11(  22,  250) $|$  729 $|$ 1061 $|$ 20000 & 7.1e-13 $|$ 1.0e-12 $|$ 7.8e-13 $|$ 4.2e-12 &06$|$13$|$18$|$7:41  & 7.7e-06 \\ \hline
				chr20b &	401 &    8(  16,  250) $|$  671 $|$ 1101 $|$ 20000 & 9.2e-13 $|$ 1.0e-12 $|$ 8.8e-13 $|$ 2.7e-12 &06$|$12$|$18$|$7:43  & 7.7e-06 \\ \hline
				chr20c &	401 &   12(  25,  230) $|$  743 $|$ 1161 $|$ 20000 & 7.1e-13 $|$ 9.8e-13 $|$ 7.2e-13 $|$ 4.4e-12 &06$|$13$|$19$|$7:35  & 7.4e-06 \\ \hline
				chr22a &	485 &   12(  28,  314) $|$  811 $|$ 1421 $|$ 20000 & 8.5e-13 $|$ 1.0e-12 $|$ 6.6e-13 $|$ 2.8e-12 &11$|$22$|$34$|$11:02 & 3.8e-06 \\ \hline
				chr22b &	485 &   10(  22,  194) $|$  743 $|$ 1261 $|$ 20000 & 1.0e-12 $|$ 1.0e-12 $|$ 8.7e-13 $|$ 2.2e-12 &08$|$19$|$30$|$11:02 & 3.8e-06 \\ \hline
				chr25a &	626 &   10(  19,  165) $|$  871 $|$ 1281 $|$ 20000 & 9.3e-13 $|$ 1.0e-12 $|$ 8.2e-13 $|$ 2.0e-12 &10$|$37$|$50$|$17:43 & 2.8e-06 \\ \hline
		\end{tabular}}
	\end{table}
\end{footnotesize}

\begin{footnotesize}
	\begin{table}[H]
		\caption{Numerical results on nug-data for QAP problems. In this table, $ y = 5\times 10^5 $.}
		\label{tab-nug}
		\resizebox{\textwidth}{!}{
			\begin{tabular}{|c|c|c|c|c|c|}\hline
				\multicolumn{2}{|c|}{}	& iteration & KKT residual & time & sc\\ \hline
				problem & $n$ & alm $|$ apg $|$ admm $|$ dykstra & alm $|$ apg $|$ admm $|$ dykstra & alm $|$ apg $|$ admm $|$ dykstra & alm\\ \hline
				nug22 & 	485 &    9(  18,  109) $|$  507 $|$ 1101 $|$ 20000 & 5.6e-13 $|$ 9.6e-13 $|$ 7.6e-13 $|$ 7.1e-12 &05$|$12$|$25$|$11:06   & 1.7e-05 \\ \hline
				nug24 & 	577 &   11(  25,  370) $|$ 1091 $|$ 1361 $|$ 20000 & 9.4e-13 $|$ 1.0e-12 $|$ 9.9e-13 $|$ 2.2e-12 &17$|$42$|$44$|$14:56   & 1.6e-06 \\ \hline
				nug25 & 	626 &   10(  21,  166) $|$ 1099 $|$ 1361 $|$ 20000 & 9.9e-13 $|$ 1.0e-12 $|$ 1.0e-12 $|$ 2.1e-12 &10$|$45$|$49$|$16:31   & 1.6e-06 \\ \hline
				nug27 & 	730 &   10(  20,  130) $|$  606 $|$ 1221 $|$ 20000 & 9.4e-13 $|$ 9.6e-13 $|$ 6.2e-13 $|$ 6.9e-12 &11$|$33$|$1:03$|$23:38 & 1.1e-05 \\ \hline
				nug30 & 	901 &    8(  19,  103) $|$  455 $|$  841 $|$ 14594 & 8.1e-13 $|$ 8.4e-13 $|$ 8.0e-13 $|$ 1.0e-12 &17$|$32$|$1:06$|$26:39 & 1.0e-04 \\ \hline
		\end{tabular}}
	\end{table}
\end{footnotesize}

\begin{footnotesize}
	\begin{table}[H]
		\caption{Numerical results on tai-data for QAP problems. In this table, $ y = 7\times  10^7 $.}
		\label{tab-tai}
		\resizebox{\textwidth}{!}{
			\begin{tabular}{|c|c|c|c|c|c|}\hline
				\multicolumn{2}{|c|}{}	& Iteration & KKT residual & time & sc\\ \hline
				problem &	n & alm $|$ apg $|$ admm $|$ dykstra & alm $|$ apg $|$ admm $|$ dykstra & alm $|$ apg $|$ admm $|$ dykstra & alm\\ \hline
				tai20b &	401 &    4(  14,   92) $|$  337 $|$  441 $|$ 1496 & 4.6e-13 $|$ 7.1e-13 $|$ 3.3e-13 $|$ 1.0e-12 &03$|$05$|$07$|$31         & 1.9e-03  \\ \hline
				tai25b &	626 &    3(   8,  104) $|$  370 $|$  401 $|$ 1017 & 4.2e-13 $|$ 2.8e-13 $|$ 5.0e-13 $|$ 1.0e-12 &06$|$13$|$16$|$49         & 4.7e-03  \\ \hline
				tai30b &	901 &    2(   7,  130) $|$  360 $|$  481 $|$ 1195 & 2.1e-13 $|$ 1.0e-13 $|$ 2.1e-13 $|$ 1.0e-12 &13$|$27$|$40$|$2:04       & 3.5e-03  \\ \hline
				tai35b &	1226 &    2(   5,  154) $|$  330 $|$  541 $|$ 1526 & 8.0e-13 $|$ 9.2e-13 $|$ 3.9e-13 $|$ 1.0e-12 &28$|$49$|$1:34$|$5:25    & 1.6e-03  \\ \hline
				tai40b &	1601 &    3(  10,  118) $|$  370 $|$  481 $|$ 1578 & 4.5e-13 $|$ 7.0e-13 $|$ 6.7e-13 $|$ 1.0e-12 &53$|$1:45$|$2:27$|$11:46 & 3.0e-03  \\ \hline
		\end{tabular}}
	\end{table}
\end{footnotesize}

\subsection{Determination of degeneracy status}

We end this section by discussing  how to determine the degeneracy status of a feasible solution $ X $ of the DNN projection problem, given the eigenvalue decomposition  $ X = PDP^T $. To this end, we may use equation~\eqref{prop:licq:eqiv} in Proposition \ref{prop:LICQ} to check whether the linear system generated by $P_\alpha^T (H_{\overline{\cal E}}+ H_{\overline{\cal E}}^T)P$ has a zero null space. In particular, if the necessary condition for constraint nondegeneracy in Proposition \ref{prop:LICQ} fails to hold, one can immediately conclude that $X$ is degenerate. However, if the necessary condition holds, then one needs to proceed to check whether the coefficient matrix generated by$P_\alpha^T (H_{\overline{\cal E}}+ H_{\overline{\cal E}}^T)P$ has full column rank. But note that since the size of the coefficient matrix is $(n|\alpha|) \times |\overline{\cal E}|$, which can be huge when $n|\alpha|$ and $|\overline{\cal E}|$ are large, it is generally expensive to numerically check the degeneracy status of $X$ in the latter case.

Table \ref{eq-test-degenearcy} presents the degeneracy status of the computed solutions for some tested instances in Table \ref{tab-Zero}--Table \ref{tab-tai}. From the table, we can see that the sizes of the corresponding linear systems are usually huge and checking whether the coefficient matrix has full column rank could be very expensive numerically.

By analysing the computational results presented in Table 1--Table 9 and the degeneracy status in Table \ref{eq-test-degenearcy}, we can observe that it is indeed more challenging to solve degenerate DNN projection problems than non-degenerate problems.

\begin{footnotesize}
	\begin{table}[htb!]
		\centering
		\caption{Degeneracy status for some tested instances in Table \ref{tab-Zero}--Table \ref{tab-tai}.}
		\label{eq-test-degenearcy}
		\resizebox{0.90\textwidth}{!}{
			\begin{tabular}{|l|c|c|c|c|c|l|}\hline
				Problems & $ n $ & $ |\alpha| $ & $\frac{1}{2}(n - |\alpha|)(n - |\alpha| +1)$ & $ |{\mathcal{E}}| $ & {\tt size of lin. sys.} & {\tt Degeneracy}  \\ \hline
				Table 1 (Zero) & 400 & 0 & 80200 & 0 & $ 0\times 80200 $ & {\tt Yes} \\ \hline
				Table 2 (Hankel) & 400 & 105 & 43660 & 37666 & $ 42000\times 42534 $ & {\tt Yes} \\ \hline
				Table 3 (low rank sparse) & 400 & 184 & 23436 & 8093  & $ 73600 \times 72107 $ & {\tt Yes} \\ \hline
				Table 4 (Toeplitz) & 400 & 357 & 946   & 11643 & $ 142800 \times 68557 $ & {\tt Out of mem.} \\ \hline
				Table 5 (bqp250-2) & 501 & 23 & 114481 & 41889 & $ 11523 \times 83862 $ & {\tt Yes} \\ \hline
				Table 6 (bur26a) & 677 & 626 & 1326 & 212603 & $ 423802\times 16900 $ & {\tt No} \\ \hline
				Table 7 (chr20a) & 401 & 362 & 780 & 73001  & $ 145162\times 7600 $ & {\tt No} \\ \hline
				Table 8 (nug22)  & 485 & 442 & 946 & 107691 & $ 214370\times 10164 $ & {\tt No} \\ \hline
				Table 9 (tai20b) & 401 & 362 & 780 & 73001 & $ 145162\times 7600 $ & {\tt No} \\ \hline
		\end{tabular}}
	\end{table}
\end{footnotesize}

\section{Conclusions}\label{sec:conclusions}
In this paper, we have employed the augmented Lagrangian method (ALM) to compute the projection onto the doubly nonnegative (DNN) cone. The ALM solves a sequence of well-conditioned nonsmooth equations instead of directly dealing with the possibly singular Karush-Kuhn-Tucker system. Under the dual quadratic growth condition and proper stopping criteria for the subproblems, the proposed algorithm is shown to converge asymptotically superlinearly.  Extensive numerical results demonstrate that our proposed ALM is more efficient and robust than the accelerated proximal gradient method, the alternative direction method of multiplier and Dykstra's algorithm. With the important role played by completely positive cone in modeling nonconvex quadratic optimization problems in various applications, we believe that our solver for computing the projection onto the DNN cone  will serve as a fundamental toolbox to approximately solve computationally intractable completely positive or copositive cone programming problems in the future.

%
%
%

\section*{Acknowledgments.}
The third author is supported in part by the Hong Kong Research Grant Council grant PolyU 153014/18P, and the fourth author is supported in part by the Ministry of Education, Singapore, under its Academic Research Fund Tier 3 grant call (MOE-2019-T3-1-010). We also thank the referees for their helpful suggestions.

\bibliographystyle{informs2014} 



\end{document}